\documentclass[11pt]{article}
\usepackage{pxfonts}
\usepackage{yfonts}
\usepackage{dsfont}
\usepackage{marvosym}
\usepackage{graphicx}
\usepackage{relsize}

\def\bel{\begin{equation}\label}
\def\eeq{\end{equation}}
\def\ds{\displaystyle}
\def\endproof{\hphantom{MM}
\hfill\llap{$\square$}\goodbreak}

\def\mt{\longrightarrow}
\def\v{\vskip 1em}

\def\ve{\varepsilon}
\def\R{\mathbb R}
\def\Z{\mathbb Z}
\def\C{\mathfrak{C}}
\def\Cx{\mathds C}

\def\Re{{\bf Re}}

\def\S{{\bf S}}

\def\Q{{\bf Q}}

\def\J{{\bf J}}

\def\L{{\bf L}}
\def\U{{\bf U}}
\def\V{{\bf V}}

\def\i{{\bf i}}
\def\Tilde{\widetilde}
\def\Hat{\widehat}

\def\bar{\overline}

\def\I{{\bf I}}

\def\M{{\bf M}}
\def\G{{\bf G}}

\def\Cup{{\bigcup}}
\def\Cap{{\bigcap}}

\def\alpha{\alphaup}
\def\beta{\betaup}
\def\gamma{\gammaup}
\def\delta{\deltaup}
\def\xi{{\xiup}}
\def\eta{{\etaup}}
\def\tau{{\tauup}}
\def\rho{{\rhoup}}
\def\phi{{\phiup}}
\def\psi{{\psiup}}
\def\lambda{{\lambdaup}}
\def\omega{\omegaup}
\def\varphi{{\varphiup}}
\def\gamma{{\gammaup}}

\def\t{{\bf t}}
\def\s{{\bf s}}
\def\r{{\bf r}}
\def\h{{\bf h}}

\newtheorem{thm}{Theorem}[section]

\newtheorem{prop}{Proposition}[section]
\newtheorem{remark}{Remark}[section]

\begin{document}
 \[\begin{array}{cc}\hbox{\LARGE{\bf Hardy-Littlewood-Sobolev inequality on product spaces}}
 \end{array}\]
 
 \[\hbox{Zipeng Wang}\]

 \begin{abstract}
We study a family of fractional integral operators defined on an homogeneous space with   a {\it rectangle doubling} - measure.  As a result, we give an extension of the classical  Hardy-Littlewood-Sobolev theorem to a multi-parameter setting.
\end{abstract}

\section{Introduction}
\setcounter{equation}{0}
In 1928, Hardy and Littlewood \cite{Hardy-Littlewood}  have first established a regularity theorem of fractional integration on $\R$. Ten years later, Sobolev \cite{Sobolev} extended this result into   higher dimensions. Today, it is well known as the Hardy-Littlewood-Sobolev inequality. 

Let $0<\alpha<n$. A fractional integral operator $I_\alpha$ is initially defined by
\bel{I_a Lebesgue}
\Big(I_\alpha f\Big)(x)~\doteq~\int_{\R^n} f(y) \left({1\over |x-y|}\right)^{n-\alpha} dy.
\eeq
$\diamond$ {\small Throughout, we regard $\C$ as a generic constant depending on its subindices}.

{\bf Hardy-Littlewood-Sobolev theorem}~~~~~{\it Let $I_\alpha$ defined in (\ref{I_a Lebesgue}) for $0<\alpha<n$. We have
\bel{norm ineq}
\begin{array}{cc}\ds
\left\| I_\alpha f\right\|_{\L^q(\R^n)}~\leq~\C_{p~q}~\left\|f\right\|_{\L^p(\R^n)},\qquad 1<p<q<\infty
\\\\ \ds
\hbox{if and only if}\qquad
{\alpha\over n}~=~{1\over p}-{1\over q}.
\end{array}
\eeq}

This classical result became fundamentally important in the application
 of many elliptic partial differential equations. Moreover, it
  has been extended to the space of homogeneous type  for certain sub-elliptic problems whose corresponding pseudo differential operators have a non-isotropic metric. For example, see the papers by Franchi and Lanconelli \cite{Franchi-Lanconelli}, Nagel, Stein and Wainger \cite{Nagel-Stein-Wainger}, Nagel \cite{Nagel} and Franchi and Serapioni \cite{Franchi-Serapioni}.

Observe that  
$\left({1\over |x-y|}\right)^{n-\alpha}=\left({1\over |x-y|^n}\right)^{1-{\alpha\over n}}$
where $|x-y|^n$ can be interpreted as the volume of the smallest ball  ( or cube ) centered on  $x$ containing $y$.
Let $\delta>0$ and 
$Q(x,\delta)\doteq \Big\{y\in\R^n~\colon~|x-y|< \delta\Big\}$.

Consider
\bel{Volume}
V(x,y)~\doteq~\inf_{\delta>0}~\Bigg\{ \mu\Big\{Q(x,\delta)\Big\}~\colon~y\in Q(x,\delta)\Bigg\}
\eeq
for which $\mu$ is an absolutely continuous measure.

We define 
\bel{I_alpha redefine}
\Big(I_\alpha f\Big)(x)~\doteq~\int_{\R^n} f(y) \left({1\over V(x,y)}\right)^{1-{\alpha\over n}} d\mu(y).
\eeq
Let $Q$ denote a cube in $\R^n$ and $2Q$ concentered with $Q$ but is doubled on the side length: $|2Q|^{1\over n}=2|Q|^{1\over n}$. A measure $\mu$ is called {\it doubling} if $d\mu(x)=\omega(x)dx$ satisfying
\bel{Doubling}
\int_{2Q} \omega(x)dx~\leq~2^{\gamma n}~ \int_Q \omega(x)dx
\eeq
for some $\gamma>0$ and every $Q\subset\R^n$. If $\mu$ is {\it doubling}, it is well known that $\mu$ is {\it reverse doubling}: 
\bel{Reverse doubling}
\int_{Q} \omega(x)dx~\leq~ 2^{-\eta n}~\int_{2Q} \omega(x)dx
\eeq
for some $\eta=\eta(\gamma)>0$ and every $Q\subset\R^n$.

{\bf Theorem One}~~~~~{\it Let $I_\alpha$ defined $w.r.t ~\mu$ in (\ref{Volume})-(\ref{I_alpha redefine}) for $0<\alpha<n$. Suppose that $\mu$ is doubling. We have
\bel{norm ineq mu}
\begin{array}{cc}\ds
\left\| I_\alpha f\right\|_{\L^q\left(\R^n, d\mu\right)}~\leq~\C_{p~q~\mu}~\left\|f\right\|_{\L^p\left(\R^n, d\mu\right)},\qquad 1<p<q<\infty
\\\\ \ds
\hbox{if and only if}\qquad
{\alpha\over n}~=~{1\over p}-{1\over q}.
\end{array}
\eeq}
\begin{remark}
The  doubling-condition (\ref{Doubling}) for $\mu$ is essential in {\bf Theorem One}. 
For instance if $\mu$ is {\it reverse doubling} satisfying (\ref{Reverse doubling}), it is allowed to be vanished on a non-zero measure set. In this case, by definition of $I_\alpha$ in (\ref{Volume})-(\ref{I_alpha redefine}), we can easily find some $f\in\L^p(\R^n)$ such that $I_\alpha f(x)$ is unbounded for every $x\in\R^n$.
\end{remark}
In practice of solving some sub-elliptic equations, we need to construct the  parametrix as a composition of two or more pseudo differential operators having different homogeneities. This situation leads us to study the operators that commute with a multi-parameter family of dilations.  During the past several decades, a number of pioneering results have been accomplished, for example  by Cordoba and Fefferman \cite{Cordoba-Fefferman}, Fefferman and Stein \cite{R.Fefferman-Stein}, Journ\'{e} \cite{Journe'}, Pipher \cite{Pipher}, Fefferman \cite{R.Fefferman} and M\"{u}ller, Ricci and Stein \cite{M.R.S}.
The area of fractional integration remains largely open. Some recent works refer to Sawyer and Wang \cite{Sawyer-Wang} and Wang \cite{Wang}.

{\bf Theorem One} can be proved by using Hedberg's method \cite{Hedberg}. A discussion is given at {\bf 4.2},  Chapter VIII of the book by Stein \cite{Stein}. In the present paper, we extend   {\bf Theorem One} to the multi-parameter setting. This answers a question proposed by professor Elias M. Stein during an informal meeting with the author in Summer 2018.

\section{Statement of the main result}
\setcounter{equation}{0}
Let $\delta$  denote an $n$-tuple $(\delta_1,\delta_2,\ldots, \delta_n)$ for $0<\delta_i<\infty, ~i=1,2,\ldots,n$ and  
\bel{Q(x,delta)}
\begin{array}{cc}\ds
\Q(x,\delta)~\doteq~\bigotimes_{i=1}^n \left\{y_i\in\R~\colon~|x_i-y_i|< \delta_i\right\}.
\end{array}
\eeq
Consider
\bel{V}
\V(x,y)~\doteq~\inf_\delta~\Bigg\{ \mu\Big\{\Q(x,\delta)\Big\}~\colon~y\in\Q(x,\delta)\Bigg\}
\eeq
which  is  the volume $w.r.t$  $\mu$ of the smallest rectangle centered on $x$ that contains $y$.
 
For $0<\alpha<n$, we define $\I_\alpha$ $w.r.t$  $\mu$  by
\bel{I_alpha}
\Big(\I_\alpha f\Big)(x)~=~\int_{\R^n} f(y) \left({1\over \V(x,y)}\right)^{1-{\alpha\over n}} d\mu(y).
\eeq
A measure $\mu$ is {\it rectangle doubling} if it is {\it doubling} on every coordinate subspace. Namely, similar to (\ref{Doubling})-(\ref{Reverse doubling}), we have
\bel{rectangle doubling}
\begin{array}{lr}\ds
\qquad\int_{2\Q_i} \omega\left(x_i,x_i'\right)dx_i~\leq~2^{\gamma n} ~\int_{\Q_i} \omega\left(x_i,x_i'\right)dx_i,
\\\\ \ds
\Longrightarrow\int_{\Q_i} \omega\left(x_i,x_i'\right)dx_i~\leq~2^{-\eta n} ~\int_{2\Q_i} \omega\left(x_i,x_i'\right)dx_i
\end{array}
\qquad x_i'\in\R^{n-1}
\eeq
for every interval $\Q_i\subset\R$, $i=1,2,\ldots,n$ and some $\gamma>0, \eta=\eta(\gamma)>0$.

{\bf Theorem Two}~~
{\it Let $\I_\alpha$  defined $w.r.t ~\mu$ in (\ref{Q(x,delta)})-(\ref{I_alpha}) for $0<\alpha<n$. Suppose that $\mu$ is $rectangle~doubling$. We have
\bel{RESULT}
\begin{array}{cc}\ds
\left\|\I_\alpha f\right\|_{\L^q\left(\R^n, d\mu\right)}~\leq~\C_{p~q~\mu}~\left\|f\right\|_{\L^p\left(\R^n,  d\mu\right)},\qquad 1<p<q<\infty
\\\\ \ds
\hbox{if and only if}\qquad
{\alpha\over n}~=~{1\over p}-{1\over q}.
\end{array}
\eeq}

{\bf Sketch of proof} : The homogeneity condition ${\alpha\over n}={1\over p}-{1\over q}$ can be verified by a standard exercise of changing dilations. 
In  section 3, we introduce a new framework where $\R^n$ is decomposed into an infinitely many {\it dyadic cones}. The consisting partial operators defined on these cones are essentially  one-parameter fractional integral operators satisfying the desired regularity. 
Furthermore, they enjoy a certain property of almost orthogonality.

Our analysis is developed in the same spirt of 
 Hedberg \cite{Hedberg}. However, instead of using Hardy-Littlewood maximal operator,  we have
 \bel{M_alpha}
\Big(\M_\alpha f\Big)(x)~\doteq~\sup_\delta~  \mu\Big\{\Q(x,\delta)\Big\}^{{\alpha\over n}-1}\int_{\Q(x,\delta)} \left|f(y)\right|d\mu(y)
\eeq
playing the role.
Section 4 is devoted to  some preliminary estimates. In particular, we obtain the next result by applying a multi-parameter Cals\'{o}n embedding theorem, recently proved by Tanaka and Yabuta \cite{Tanaka-Yabuta}.
\begin{thm}
Let $\M_\alpha$ defined in (\ref{M_alpha}) for $0<\alpha<n$. We have
\bel{M regularity}
\begin{array}{cc}\ds
\left\|\M_\alpha f\right\|_{\L^q\left(\R^n, d\mu\right)}~\leq~\C_{p~q~\mu}~\left\| f\right\|_{\L^p\left(\R^n,d\mu\right)},\qquad 1<p<q<\infty
\\\\ \ds
\hbox{if}\qquad {\alpha\over n}~=~{1\over p}-{1\over q}.
\end{array}
\eeq
\end{thm}
In section 5,  we prove a crucial lemma of almost orthogonality.

\section{Cone decomposition}
\setcounter{equation}{0}
Let $\t$  denote an $n$-tuple $(2^{-t_1},2^{-t_2},\ldots,2^{-t_n})$  for  $t_i\in\Z,~i=1,2,\ldots,n$.  
We define
\bel{Partial}
\Big(\Delta_\t \I_\alpha f\Big)(x)~\doteq~\int_{\Gamma_\t(x)} f(y)\left({1\over \V(x,y)}\right)^{1-{\alpha\over n}} d\mu(y)
\eeq
where
\bel{Cone}
\begin{array}{lr}\ds
\Gamma_\t(x)~\doteq~\Cup_{\jmath\in\Z}~\bigotimes_{i=1}^n\left\{y_i\in\R\colon~2^{\jmath-t_i}\leq |x_i-y_i|< 2^{\jmath+1-t_i}\right\}.
\end{array}
\eeq
Observe that $\Gamma_\t(x)$ in (\ref{Cone}) is a collection of rectangles having  the same eccentricity  $w.r.t ~\t$, whose diameters are comparable to their distances away from $x$, in the sprit of Witney.  Geometrically, it is interpreted as a {\it dyadic cone} with vertex on $x$.

In particular, suppose  $\Hat{\t}$ is another $n$-tuple $(2^{-t_1-\ell}, 2^{-t_2-\ell},\ldots,2^{-t_n-\ell}),~\ell\in\Z$. We must have $\Gamma_\t(x)=\Gamma_{\Hat{\t}}(x)$ since the union in (\ref{Cone}) takes all $\jmath\in\Z$. 
\begin{remark}
Without lose of the generality, we assume  $t_i,i=1,2,\ldots,n$ to be nonnegative integers and
$t_\nu\doteq\min\{t_i\colon~i=1,2,\ldots,n\}=0$.
In the case of $t_1=t_2=\cdots=t_n=0$, we write  $\t=o$.
\end{remark}
\begin{figure}[h]
\centering
\includegraphics[scale=0.55]{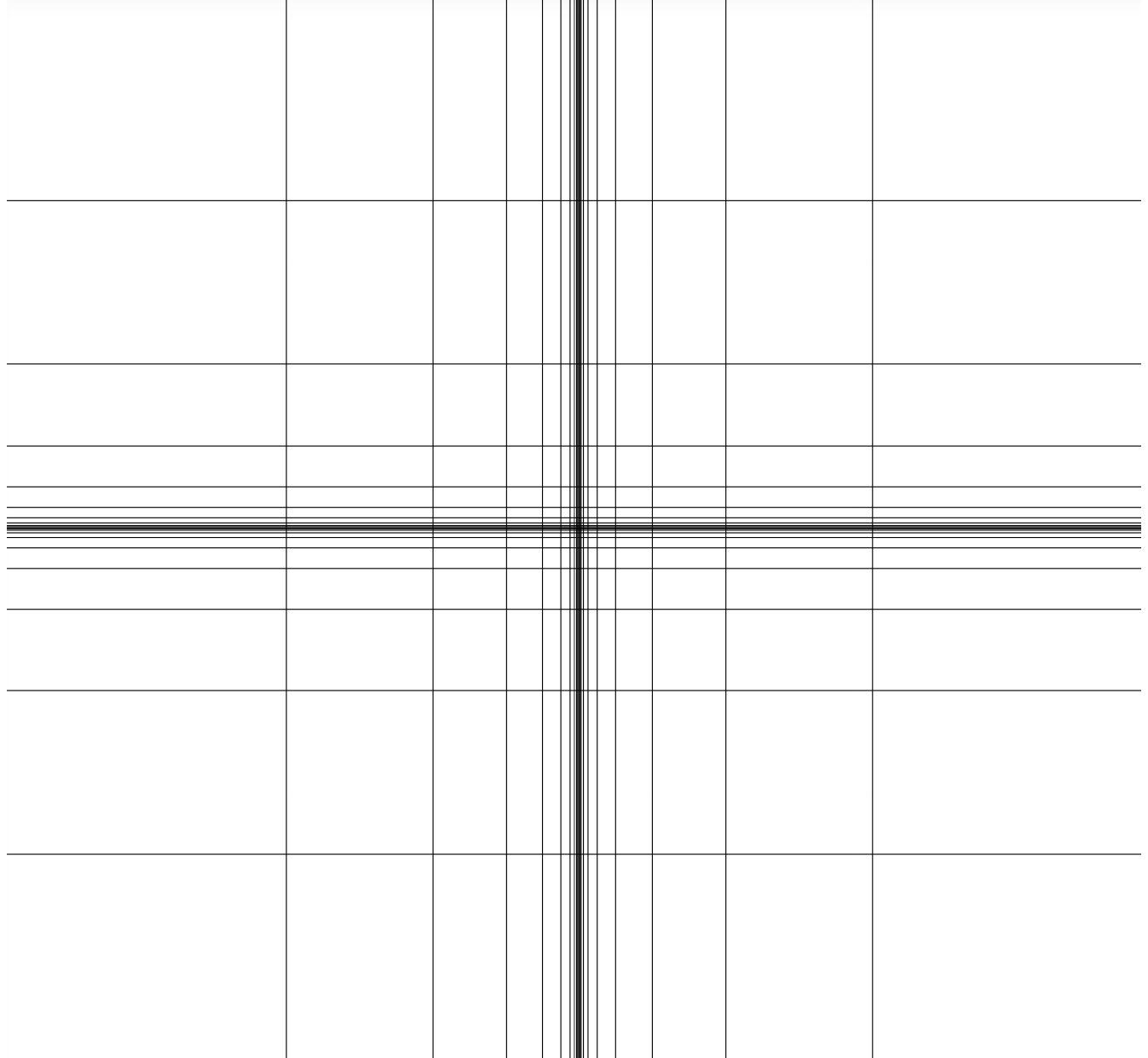}
\caption{The collection of all dyadic rectangles in a $2$-parameter setting.}
\end{figure}

$\diamond$ For brevity, we abbreviate $\L^p\left(\mu\right)=\L^p\left(\R^n, d\mu\right)$.

Given $f\ge0\in\L^p(\mu)$,  define
\bel{theta_t}
\vartheta_\t(x)~=~  {1\over  \left\| f\right\|^p_{\L^p\left(\mu\right)} }         \int_{\Gamma_\t(x)} \Big(f(y)\Big)^p d\mu(y).
\eeq
Let ${\alpha\over n}={1\over p}-{1\over q},~1<p<q<\infty$.  We  claim   
\bel{Regularity est}
\begin{array}{cc}\ds
\Big(\Delta_\t \I_\alpha f\Big)(x)
~\leq~\C_{p~q~\mu}~\Big(\vartheta_\t(x)\Big)^{{1\over p}\left(1-{r\over q}\right)}\Big(\M_\beta f\Big)^{r\over q}(x)\left\| f\right\|_{\L^p\left(\mu\right)}^{1-{r\over q}},
\\\\ \ds
{\beta\over n}~=~{1\over p}-{1\over r},\qquad p<r<q
\end{array}
\eeq
where $\M_\beta$ is defined in (\ref{M_alpha}).

By definition of $\V(x,y)$ in (\ref{V}), we have
\bel{V norm compara}
\begin{array}{rl}\ds
\mu\left\{\bigotimes_{i=1}^n|x_i-y_i|\leq 2^{\jmath-t_i}\right\} ~\leq~\V(x,y)~<~\mu\left\{\bigotimes_{i=1}^n|x_i-y_i|\leq 2^{\jmath+1-t_i}\right\}
\end{array}
\eeq
whenever
\bel{Gamma_t^j}
y~\in~ \Gamma_\t^\jmath(x)~\doteq~\bigotimes_{i=1}^n\left\{y_i\in\R\colon~2^{\jmath-t_i}\leq |x_i-y_i|< 2^{\jmath+1-t_i}\right\}.
\eeq
Define $\tau(\t,x)\in\R$ implicitly by
\bel{tau(t,x)}
\begin{array}{lr}\ds
{\Big(\vartheta_\t(x)\Big)^{1\over p}\over \Big(\M_\beta f\Big)(x)}\left\| f\right\|_{\L^p\left(\mu\right)}~\doteq~\mu\left\{\bigotimes_{i=1}^n |x_i-y_i|\leq2^{\tau(\t,x)-t_i}\right\}^{1\over r}.
\end{array}
\eeq
Suppose that $\mu$ is {\it rectangle doubling} satisfying (\ref{rectangle doubling}). 
Consider $-\infty<\jmath\leq\tau(\t,x)$. We have
\bel{Region-I Est}
\begin{array}{lr}\ds
 \int_{\Gamma_\t^\jmath(x)} f(y)\left({1\over \V(x,y)}\right)^{1-{\alpha\over n}} d\mu(y)
~\leq~\mu\left\{\bigotimes_{i=1}^n |x_i-y_i|\leq 2^{\jmath-t_i}\right\}^{{\alpha\over n}-1} \int_{\Gamma_\t^\jmath(x)} f(y)d\mu(y)
 \\\\ \ds
~\leq~
 \mu\left\{\bigotimes_{i=1}^n |x_i-y_i|\leq 2^{\jmath-t_i}\right\}^{\alpha-\beta\over n}~{2^{\gamma n\left(1-{\beta\over n}\right)}\over \ds \mu\left\{\bigotimes_{i=1}^n |x_i-y_i|\leq 2^{\jmath+1-t_i}\right\}^{1-{\beta\over n}}} \int_{\Gamma_\t^\jmath(x)} f(y)d\mu(y)
\\\\ \ds
~\leq~2^{\gamma n\left(1-{\beta\over n}\right)}\mu\left\{\bigotimes_{i=1}^n |x_i-y_i|\leq2^{\jmath-t_i}\right\}^{\alpha-\beta\over n}\Big(\M_\beta f\Big)(x) 
\\\\ \ds
~\leq~2^{-\eta(\alpha-\beta)(\tau(\t,x)-\jmath)} 2^{\gamma n\left(1-{\beta\over n}\right)} \mu\left\{\bigotimes_{i=1}^n |x_i-y_i|\leq2^{\tau(\t,x)-t_i}\right\}^{\alpha-\beta\over n}\Big(\M_\beta f\Big)(x).\qquad(\alpha-\beta>0)
\end{array}
\eeq
Consider $\tau(\t,x)\leq\jmath<\infty$. By applying H\"{o}lder inequality, we have
\bel{Region-II Holder}
\begin{array}{lr}\ds
\int_{ \Gamma_\t^\jmath(x)}  f(y)\left({1\over \V(x,y)}\right)^{1-{\alpha\over n}} d\mu(y)
~\leq~\Big(\vartheta_\t(x)\Big)^{1\over p} \left\| f\right\|_{\L^p\left(\mu\right)}
\left\{\int_{ \Gamma_\t^\jmath(x)  }   \left({1\over \V(x,y)}\right)^{\left(1-{\alpha\over n}\right)\left({p\over p-1}\right)} d\mu(y)\right\}^{p-1\over p}
\end{array}
\eeq
where
\bel{Region-II Est}
\begin{array}{lr}\ds
\int_{ \Gamma_\t^\jmath(x)  }  \left({1\over \V(x,y)}\right)^{\left(1-{\alpha\over n}\right)\left({p\over p-1}\right)} d\mu(y)
~\leq~  \mu\left\{\bigotimes_{i=1}^n |x_i-y_i|\leq 2^{\jmath-t_i}\right\}^{\left({\alpha\over n}-1\right)\left({p\over p-1}\right)} \mu\left\{\bigotimes_{i=1}^n|x_i-y_i|\leq2^{\jmath+1-t_i}\right\}   
\\\\ \ds~~~~~~~~~
~\leq~2^{\gamma n}\mu\left\{\bigotimes_{i=1}^n|x_i-y_i|\leq2^{\jmath-t_i}\right\}^{\left({\alpha\over n}-{1\over p}\right)\left({p\over p-1}\right)}
\\\\ \ds~~~~~~~~~
~\leq~2^{-\eta(\jmath-\tau(\t,x))\left({ n\over q}\right)\left({p\over p-1}\right)} 2^{\gamma n}\mu\left\{\bigotimes_{i=1}^n|x_i-y_i|\leq2^{\tau(\t,x)-t_i}\right\}^{\left({{\alpha\over n}-{1\over p}}\right)\left({p\over p-1}\right)}.\qquad\hbox{\small{$\left({\alpha\over n}-{1\over p}=-{1\over q}<0\right)$}}
\end{array}
\eeq
By bringing (\ref{tau(t,x)}) to  (\ref{Region-I Est}), we find
\bel{Regularity Est I}
\begin{array}{lr}\ds
\int_{\Gamma_\t^\jmath(x)}f(y)  \left({1\over \V(x,y)}\right)^{1-{\alpha\over n}} d\mu(y)
~\leq~2^{-\eta(\alpha-\beta)(\tau(\t,x)-\jmath)} 2^{\gamma n\left(1-{\beta\over n}\right)} \mu\left\{\bigotimes_{i=1}^n |x_i-y_i|\leq2^{\tau(\t,x)-t_i}\right\}^{\alpha-\beta\over n} \Big(\M_\beta f\Big)(x)
\\\\ \ds~~~~~~~~~~~~~~~~~~~~~~~~~~~~~~~~~~~~~~~~~~~~~~~
~=~2^{-\eta(\alpha-\beta)(\tau(\t,x)-\jmath)} 2^{\gamma n\left(1-{\beta\over n}\right)}\left\{ {\Big(\M_\beta f\Big)(x)\over\Big(\vartheta_\t(x)\Big)^{1\over p}\left\| f\right\|_{\L^p\left(\mu\right)}}\right\}^{{r\over q}-1}\Big(\M_\beta f\Big)(x)
\\\\ \ds~~~~~~~~~~~~~~~~~~~~~~~~~~~~~~~~~~~~~~~~~~~~~~~
~=~2^{-\eta(\alpha-\beta)(\tau(\t,x)-\jmath)} 2^{\gamma n\left(1-{\beta\over n}\right)}\Big(\vartheta_\t(x)\Big)^{{1\over p}\left(1-{r\over q}\right)}\Big(\M_\beta f\Big)^{r\over q}(x)\left\| f\right\|_{\L^p\left(\mu\right)}^{1-{r\over q}}.
\end{array}
\eeq
By bringing (\ref{tau(t,x)}) to (\ref{Region-II Holder})-(\ref{Region-II Est}) 
, we find
\bel{Regularity Est II}
\begin{array}{lr}\ds
\int_{\Gamma_\t^\jmath(x)}f(y)  \left({1\over \V(x,y)}\right)^{1-{\alpha\over n}} d\mu(y)
\\\\ \ds
~\leq~2^{-\eta(\jmath-\tau(\t,x))\left({ n\over q}\right)} 2^{\gamma n\left({p-1\over p}\right)}\mu\left\{ \bigotimes_{i=1}^n |x_i-y_i|\leq2^{\tau(\t,x)-t_i}\right\}^{{\alpha\over n}-{1\over p}}\Big(\vartheta_\t(x)\Big)^{1\over p} \left\| f\right\|_{\L^p\left(\mu\right)}
\\\\ \ds
~=~2^{-\eta(\jmath-\tau(\t,x))\left({n\over q}\right)} 2^{\gamma n\left({p-1\over p}\right)}\left\{ {\Big(\M_\beta f\Big)(x)\over\Big(\vartheta_\t(x)\Big)^{1\over p}\left\| f\right\|_{\L^p\left(\mu\right)}}\right\}^{r\over q}\Big(\vartheta_\t(x)\Big)^{1\over p}\left\| f\right\|_{\L^p\left(\mu\right)}
\\\\ \ds
~=~2^{-\eta(\jmath-\tau(\t,x))\left({ n\over q}\right)}2^{\gamma n\left({p-1\over p}\right)} \Big(\vartheta_\t(x)\Big)^{{1\over p}\left(1-{r\over q}\right)}\Big(\M_\beta f\Big)^{r\over q}(x)\left\| f\right\|_{\L^p\left(\mu\right)}^{1-{r\over q}}.
\end{array}
\eeq
By putting together (\ref{Regularity Est I}) and (\ref{Regularity Est II}),  we find 
\bel{Regularity Est Rectangle}
\begin{array}{lr}\ds
\int_{\Gamma_\t^\jmath(x)} f(y)\left({1\over \V(x,y)}\right)^{1-{\alpha\over n}} d\mu(y)
\\\\ \ds
~\leq~2^{-\eta|\jmath-\tau(\t,x)|\min\left\{\alpha-\beta, {n\over q} \right\}} 2^{\gamma n}\Big(\vartheta_\t(x)\Big)^{{1\over p}\left(1-{r\over q}\right)}\Big(\M_\beta f\Big)^{r\over q}(x)\left\| f\right\|_{\L^p\left(\mu\right)}^{1-{r\over q}},\qquad \jmath\in\Z.
\end{array}
\eeq
From (\ref{Cone}), we have $\Gamma_\t(\xi)=\sum_{\jmath\in\Z} \Gamma_\t^\jmath(x)$. By  summing over all $\jmath\in\Z$ in (\ref{Regularity Est Rectangle}), we obtain (\ref{Regularity est}).

Let $\t-\h$ be another $n$-tuple as $\t$. Our main objective is to prove the following result.

{\bf Lemma of almost orthogonality:} ~~{\it Suppose that $\mu$ is $rectangle~doubling$. We have
\bel{Ortho Result}
\begin{array}{cc}\ds
\int_{\R^n}\sum_\t \Big(\Delta_\t \I_\alpha f\Big)(x) \Big(\Delta_{\t-\h}\I_\alpha f\Big)^{q-1}(x)d\mu(x)
~\leq~\C_{p~q~\mu}~\prod_{i=1}^n 2^{-\ve|h_i|}  \left\| f\right\|_{\L^p\left(\mu\right)}^{q}
\\\\ \ds
\hbox{if}\qquad {\alpha\over n}~=~{1\over p}-{1\over q},\qquad 1<p<q<\infty
 \end{array}
\eeq 
for some $\ve=\ve(p,q,\mu)>0$ and $q\in\Z$ sufficiently large.}

Let $\t-\h^m,m=1,2,\ldots,q-1$ denote for the $n$-tuples as $\t-\h$. 
We have
\bel{Sum}
\begin{array}{lr}\ds
\int_{\R^n} \Big(\I_\alpha f\Big)^q(x)d\mu(x)
~=~\sum_{\h^m,m=1,2,\ldots,q-1}\int_{\R^n}\sum_\t \Big(\Delta_\t \I_\alpha f\Big)(x) \prod_{m=1}^{q-1}\Big(\Delta_{\t-\h^m}\I_\alpha f\Big)(x) d\mu(x).
\end{array}
\eeq
By applying H\"{o}lder inequality twice, we find
\bel{Ortho Expansion}
\begin{array}{lr}\ds
\int_{\R^n}\sum_\t \Big(\Delta_\t \I_\alpha f\Big)(x) \prod_{m=1}^{q-1}\Big(\Delta_{\t-\h^m}\I_\alpha f\Big)(x) d\mu(x)
~\leq~ \int_{\R^n} \prod_{m=1}^{q-1} \left\{\sum_\t \Big(\Delta_\t \I_\alpha f\Big)(x) \Big(\Delta_{\t-\h^m}\I_\alpha f\Big)^{q-1}(x)\right\}^{1\over q-1}d\mu(x)
\\\\ \ds~~~~~~~~~~~~~~~~~~~~~~~~~~~~~~~~~~~~~~~~~~~~~~~~~~~~~~~~~~~~~~~~~~~~~~
~\leq~\prod_{m=1}^{q-1} \left\{\int_{\R^n} \sum_\t \Big(\Delta_\t \I_\alpha f\Big)(x) \Big(\Delta_{\t-\h^m}\I_\alpha f\Big)^{q-1}(x)d\mu(x)\right\}^{1\over q-1}.
\end{array}
\eeq
Let $\vartheta_\t(x)$ defined in (\ref{theta_t}). We have $\sum_\t \vartheta_\t(x)=1$.
Suppose $p\leq q-r$. From (\ref{Sum})-(\ref{Ortho Expansion}) and by using (\ref{Ortho Result}), we have 
\bel{Heuristic Est}
\begin{array}{lr}\ds
\int_{\R^n} \Big(\I_\alpha f\Big)^q(x)
~\leq~\C_{p~q~\mu}\left\| f\right\|_{\L^p\left(\mu\right)}^{q-r} \int_{\R^n}\sum_\t \Big(\vartheta_\t(x)\Big)^{{q-r\over p} }\Big(\M_\beta f\Big)^{r}(x)d\mu(x)\qquad \hbox{\small{by (\ref{Regularity est})}}
\\\\ \ds~~~~~~~~~~~~~~~~~~~~~
~\leq~\C_{p~q~\mu}\left\| f\right\|_{\L^p\left(\mu\right)}^{q-r} \int_{\R^n} \left\{\sum_\t \vartheta_{\t}(x)\right\} \Big(\M_\beta f\Big)^{r}(x)d\mu(x)
 \\\\ \ds~~~~~~~~~~~~~~~~~~~~~
 ~\leq~\C_{p~q~\mu}\left\| f\right\|_{\L^p\left(\mu\right)}^{q-r} \int_{\R^n}\Big(\M_\beta f\Big)^{r}(x)d\mu(x)
\\\\ \ds~~~~~~~~~~~~~~~~~~~~~
~\leq~\C_{p~q~\mu}\left\| f\right\|_{\L^p\left(\mu\right)}^{q}\qquad\hbox{\small{by {\bf Theorem 2.1}.}}
\end{array}
\eeq
Note that
${\alpha\over n}={1\over p}-{1\over q}={q-1\over q}-{p-1\over p}$.
The lemma holds by replacing $p,q$   with 
${q\over q-1},{p\over p-1}$ respectively. 
Because $\I_\alpha$ is  self-adjoint, we have $\I_\alpha\colon \L^p\left(\mu\right)\mt\L^q\left(\mu\right)\Longleftrightarrow
\I_\alpha\colon \L^{q\over q-1}\left(\mu\right)\mt\L^{p\over p-1}\left(\mu\right)
$.

Let
\bel{formula12}
\begin{array}{rl}\ds
{1\over p}-{1\over q}~=~{\alpha\over n}~=~{1\over p_1}-{1\over q_1}~=~{1\over p_2}-{1\over q_2},
\\\\ \ds
 1~<~p_i~<~q_i~<~\infty,\qquad i~=~1,2.
 \end{array}
\eeq
Choose $q_1,\left({p_2\over p_2-1}\right)\in\Z$  sufficiently large depending on  $p_1,q_2$ and $\mu$ respectively. By keeping the equalities hold in (\ref{formula12}), there exists a $0\leq t\leq 1$ such that 
\bel{t,pq}
{1\over p}~=~{1-t\over p_1}+{t\over p_2},\qquad{1\over q}~=~{1-t\over q_1}+{t\over q_2}.
\eeq
The  estimates in (\ref{Sum})-(\ref{Heuristic Est}) imply 
simultaneously 
\bel{Endpoints Est}
\left\|\I_\alpha f\right\|_{\L^{q_1}\left(\mu\right)}~\leq~\C_{p_1~q_1~\mu}~\left\|f\right\|_{\L^{p_1}\left(\mu\right)},\qquad \left\|\I_\alpha f\right\|_{\L^{q_2}\left(\mu\right)}~\leq~\C_{p_2~q_2~\mu}~\left\|f\right\|_{\L^{p_2}\left(\mu\right)}.
\eeq  
By applying Riesz-Thorin interpolation theorem \cite{Stein-Weiss*}, we obtain the norm inequality in (\ref{RESULT}).

\section{Some preliminary estimates}
\setcounter{equation}{0}
{\bf ( 1 )} Note that $\M_\alpha$ given in (\ref{M_alpha})  can be equivalently defined by
\bel{M_alpha equi}
\Big(\M_\alpha f\Big)(x)~=~\sup_{\Q~\ni~ x}  ~\mu\left\{\Q\right\}^{{\alpha\over n}-1}\int_{\Q} f(y)d\mu(y).
\eeq
Let $\mathfrak{D}^i,~i=1,2,\ldots,n$ denote a dyadic grid \footnote{the family of all dyadic (interval) cube: $ 2^{-k}\left\{j+[0,1)\right\}\subset\R,~~k,~j\in\Z$.} in the $i$-th one-dimensional subspace. Moreover, $\mathfrak{D}=\mathfrak{D}^1\times\mathfrak{D}^2\times\cdots\times\mathfrak{D}^n$ forms a dyadic partial grid in $\R^n$. ( it fails the nested property ) 

We write $\Q\in\mathfrak{D}$ if $\Q_i\in\mathfrak{D}^i$ for every $i=1,2,\ldots,n$.
A  dyadic  strong  fractional maximal operator  $\M_\alpha^\triangle$ is defined by
\bel{M_alpha dyadic}
\begin{array}{lr}\ds
\Big(\M_\alpha^\triangle f\Big)(x)~=~\sup_{\Q~\ni~ x\colon~\Q~\in~\mathfrak{D}}  ~\mu\left\{\Q\right\}^{{\alpha\over n}-1}\int_{\Q} f(y)d\mu(y)
\\\\ \ds~~~~~~~~~~~~~~~~
~\leq~\sum_{\Q~\ni~x\colon~\Q~\in~\mathfrak{D}} \mu\left\{\Q\right\}^{{\alpha\over n}-1}\int_{\Q} f(y)d\mu(y).
\end{array}
\eeq
By using  the Str\"{o}mberg ${1\over 3}$-trick,\footnote{~~See Sawyer and Wang \cite{Sawyer-Wang} for example and many references cited there.} we can bound $\M_\alpha$  in (\ref{M_alpha equi}) by a sum  of  finitely many dyadic strong   fractional maximal operators.  Each one  of them is defined as (\ref{M_alpha dyadic}) but  on  a shifted dyadic partial grid. 
The regarding estimate is found many times over in literatures.  
Therefore, it is suffice to study  the regularity of $\M_\alpha^\triangle$ instead.

Let ${\alpha\over n}={1\over p}-{1\over q},1<p<q<\infty$. For $f\in\L^p(\mu)$ and $g\in\L^{q\over q-1}(\mu)$, we have 
\bel{Discrete Est}
\begin{array}{lr}\ds
\int_{\R^n} \left\{\sum_{\Q~\ni~x\colon~\Q~\in~\mathfrak{D}} \mu\left\{\Q\right\}^{{\alpha\over n}-1}\int_{\Q} f(y)d\mu(y)\right\} g(x)d\mu(x)
\\\\ \ds
~=~\sum_{\Q~\in~\mathfrak{D}} \mu\left\{\Q\right\}^{{\alpha\over n}-1} \left\{\int_{\Q} f(y)d\mu(y)\right\}\left\{\int_{\Q}g(x)d\mu(x)\right\}
\\\\ \ds
~=~\sum_{\Q~\in~\mathfrak{D}} \mu\left\{\Q\right\}^{{\alpha\over n}-1} \mu\left\{\Q\right\}^{1\over p} \mu\left\{\Q\right\}^{p-1\over p} \left\{{1\over \mu\left\{\Q\right\}}\int_{\Q} f(y)d\mu(y)\right\}\mu\left\{\Q\right\}^{1\over q}\mu\left\{\Q\right\}^{q-1\over q}\left\{{1\over \mu\left\{\Q\right\}}\int_{\Q}g(x)d\mu(x)\right\}
\\\\ \ds
~=~\sum_{\Q~\in~\mathfrak{D}}  \mu\left\{\Q\right\}^{1\over p} \left\{{1\over \mu\left\{\Q\right\}}\int_{\Q} f(y)d\mu(y)\right\}\mu\left\{\Q\right\}^{q-1\over q}\left\{{1\over \mu\left\{\Q\right\}}\int_{\Q}g(x)d\mu(x)\right\}
\\\\ \ds
~\leq~\left\{ \sum_{\Q~\in~\mathfrak{D}}  \mu\left\{\Q\right\}^{r\over p} \left\{{1\over \mu\left\{\Q\right\}}\int_{\Q} f(y)d\mu(y)\right\}^r\right\}^{1\over r}\left\{\sum_{\Q~\in~\mathfrak{D}} \mu\left\{\Q\right\}^{\left({q-1\over q}\right)\left({r\over r-1}\right)}\left\{{1\over \mu\left\{\Q\right\}}\int_{\Q}g(x)d\mu(x)\right\}^{r\over r-1}\right\}^{r-1\over r}
\\ \ds~~~~~~~~~~~~~~~~~~~~~~~~~~~~~~~~~~~~~~~~~~~~~~~~~~~~~~~~~~~~~~~~~~~~~~~~~~~~~~~~~~~~~~~~~~~~
~~\hbox{\small{by H\"{o}lder inequality ( $p<r<q$ ). }}
\end{array}
\eeq
Now, we recall a multi-parameter version of Cales\'{o}n embedding theorem.
\begin{thm}
{\bf Tanaka and Yabuta, 2019~~}  Suppose that $\mu$ is $rectangle~doubling$. We have
\bel{product Carleson embedding}
\left\{ \sum_{\Q~\in~\mathfrak{D}}  \mu\left\{\Q\right\}^{r\over p} \left\{{1\over \mu\left\{\Q\right\}}\int_{\Q} f(y)d\mu(y)\right\}^r\right\}^{1\over r}~\leq~\C_{p~r~\omega}~\left\|f\right\|_{\L^p\left(\mu\right)},\qquad 1<p<r<\infty.
\eeq
\end{thm}

This result is proved by using a clever iteration argument. See the paper by  Tanaka and Yabuta \cite{ Tanaka-Yabuta}. 

By applying {\bf Theorem 4.1} to the last line of (\ref{Discrete Est}) and then taking the supremum of all $ \left\|g\right\|_{\L^{q\over q-1}(\mu)}=1$, we finish the proof of  {\bf Theorem 2.1}.

{\bf ( 2 )}  In order to prove (\ref{Ortho Result}), it is suffice to consider  $\alpha$  sufficiently close to $n$, such that
\bel{Constraint Case3}
\eta-\gamma\left(1-{\alpha\over n}\right)~>~0.
\eeq
Let ${\alpha\over n}={1\over p}-{1\over q}$. Consider 
${\alpha_i\over n}={1\over p_i}-{1\over q}, i=1,2$ 
where $p_2<p<p_1$ and $\alpha_1<\alpha<\alpha_2$.
There exists a $0<t<1$ such that ${1\over p}={1-t\over p_1}+{t\over p_2}$ and  $\alpha=(1-t)\alpha_1+t\alpha_2$. 

Set $\eta-\gamma\left(1-{\alpha_2\over n}\right)=\eta-\gamma\left(1-{1\over p_2}+{1\over q}\right)>0$ by choosing $p_2$ sufficiently close to $1$ and $q\in\Z$ sufficiently large.

Define
 $\alpha_z\doteq(1-z)\alpha_1+z\alpha_2$ and $f_z\doteq f^{p\left[{1-z\over p_1}+{z\over p_2}\right]}$ on the strip $\S\doteq\{z\in\Cx~\colon~ 0<\Re z<1\}$. Moreover, 
 
for $z=\lambda+\i v, \lambda, v\in\R$, write
${1\over p_\lambda}\doteq{1-\lambda\over p_1}+{\lambda\over p_2}$, $\alpha_\lambda\doteq\Re \alpha_z$  and $f_\lambda\doteq|f_z|$. In particular,  we have $\|f_\lambda\|_{\L^{p_\lambda}\left(\mu\right)}=\|f\|_{\L^p\left(\mu\right)}=1$. 

For every $\t$, we define
\bel{J(x,z)}
\Big(\Delta_\t\J_{\alpha~z} f\Big)(x)~\doteq~\int_{\Gamma_\t(x)} f_z(y)\left({1\over \V(x,y)}\right)^{1-{\alpha_z\over n}} d\mu(y),\qquad z\in\S.
\eeq
Let ${\beta_\lambda\over n}={1\over p_\lambda}-{1\over r_\lambda}$ for $p_\lambda<r_\lambda<q, ~0\leq\lambda\leq1$. By using (\ref{Regularity est}), we have
\bel{J(x,z), norm}
\begin{array}{lr}\ds
\left|\Big(\Delta_\t\J_{\alpha~z} f\Big)(x)\right|~\leq~\int_{\Gamma_\t(x)} f_\lambda(y)\left({1\over \V(x,y)}\right)^{1-{ \alpha_\lambda\over n}} d\mu(y)\qquad z\in\bar{\S}
\\\\ \ds~~~~~~~~~~~~~~~~~~~~~~
~\leq~\C_{p_\lambda~q~\gamma~\eta}~\Big(\vartheta_\t(x)\Big)^{{1\over p_\lambda}\left(1-{r_\lambda\over q}\right)}
\Big(\M_{\beta_\lambda} f_\lambda \Big)^{r_\lambda\over q}(x)\left\| f_\lambda\right\|_{\L^{p_\lambda}\left(\mu\right)}^{1-{r_\lambda\over q}}.
\end{array}
\eeq

For $a.e$ $x,y\in\R^n$, the integrand  in (\ref{J(x,z)}) is analytic. Together with (\ref{J(x,z), norm}), it follows that  $\left(\Delta_\t\J_{\alpha~z} f\right)(x)$ itself has a power series expansion at   every  $z\in\bar{\S}$.

Next, we consider
\bel{U_h}
\Big(\U_{\alpha~\h} f\Big)(z)~\doteq~\int_{\R^n}\sum_\t \Big(\Delta_\t\J_{\alpha~z} f\Big)(x)\Big(\Delta_{\t-\h}\J_{\alpha~z}f\Big)^{q-1}(x)d\mu(x).
\eeq
Suppose that $q\in\Z$ is sufficiently large:
$\left({q-2\over p_\lambda}\right)\left({q-r\over q}\right)\ge1$.
Recall $\vartheta_\t(x)$ defined in (\ref{theta_t}) whereas $\sum_\t \vartheta_\t(x)=1$.
By using  (\ref{J(x,z), norm}), we have
\bel{U_h bound}
\begin{array}{lr}\ds
\left|\Big(\U_{\alpha~\h} f\Big)(z)\right|~\leq~\int_{\R^n}\sum_\t \left|\Big(\Delta_\t\J_{\alpha~z} f\Big)(x)\right|\left|\Big(\Delta_{\t-\h}\J_{\alpha~z}f\Big)^{q-1}(x)\right|d\mu(x)\qquad z\in\bar{\S}
\\\\ \ds
~\leq~\C_{p_\lambda~q~\gamma~\eta}~\left\| f_\lambda\right\|_{\L^{p_\lambda}\left(\mu\right)}^{q-r} \int_{\R^n}\sum_\t \Big(\vartheta_\t(x)\Big)^{{1\over p_\lambda}\left(1-{r\over q}\right)}\Big(\vartheta_{\t-\h}(x)\Big)^{{q-1\over p_\lambda}\left(1-{r\over q}\right)}
\Big(\M_{\beta_\lambda} f_\lambda\Big)^{r}(x)d\mu(x)
\\\\ \ds
~\leq~\C_{p_\lambda~q~\gamma~\eta}~\left\| f_\lambda\right\|_{\L^{p_\lambda}\left(\mu\right)}^{q-r} \int_{\R^n} \left\{\sum_\t\Big(\vartheta_{\t-\h}(x)\Big)^{\left({q-2\over p_\lambda}\right)\left({q-r\over q}\right)}\right\} \Big(\M_{\beta_\lambda} f_\lambda\Big)^{r}(x)d\mu(x)
\\\\ \ds
~\leq~\C_{p_\lambda~q~\gamma~\eta}~\left\| f_\lambda\right\|_{\L^{p_\lambda}\left(\mu\right)}^{q-r} \int_{\R^n} \left\{\sum_\t\vartheta_{\t-\h}(x)\right\} \Big(\M_{\beta_\lambda} f_\lambda\Big)^{r}(x)d\mu(x)
 \\\\ \ds
 ~\leq~\C_{p_\lambda~q~\gamma~\eta}~\left\| f_\lambda\right\|_{\L^{p_\lambda}\left(\mu\right)}^{q-r} \int_{\R^n}\Big(\M_{\beta_\lambda} f_\lambda\Big)^{r}(x)d\mu(x)
 \\\\ \ds
~\leq~\C_{p_\lambda~q~\gamma~\eta}~\left\| f_\lambda\right\|_{\L^{p_\lambda}\left(\mu\right)}^{q}\qquad\hbox{\small{by {\bf Theorem 2.1}}}
\\\\ \ds
~=~\C_{p_\lambda~q~\gamma~\eta}.\qquad \hbox{\small{ ( $\| f_\lambda\|_{\L^{p_\lambda}}=1$ )}}
\end{array}
\eeq
Observe that $(\V_{\alpha~\h} f)(z)$ is analytic for $z\in\S$ and continuous on $\bar{\S}$.  

Suppose  that (\ref{Ortho Result})  can be obtained for ${\alpha_2\over n}={1\over p_2}-{1\over q},~1<p_2<q<\infty$ where  $\eta-\gamma\left(1-{\alpha_2\over n}\right)>0$ and $q$ sufficiently large: $\left({q-2\over p_\lambda}\right)\left({q-r\over q}\right)\ge1$. From  (\ref{U_h})-(\ref{U_h bound}), we simultaneously have
\bel{0,1,Est}
\Big(\U_{\alpha~\h} f\Big)(0+\i v)~\leq~\prod_{i=1}^n 2^{-\ve|h_i|}~\C_{p_2~q~\gamma~\eta},\qquad \Big(\U_{\alpha~\h} f\Big)(1+\i v)~\leq~\C_{p_1~q~\gamma~\eta}.
\eeq
By applying the Three-Line Lemma \cite{Stein-Weiss*}, we obtain
\bel{U_h(t)}
\begin{array}{lr}\ds
\Big(\U_{\alpha~\h} f\Big)(t+\i v)~\leq~\prod_{i=1}^n 2^{-\ve t|h_i|}~\left(\C_{p_1~q~\gamma~\eta}\right)^{1-t}\left(\C_{p_2~q~\gamma~\eta}\right)^t
\\\\ \ds~~~~~~~~~~~~~~~~~~~~~~~
~=~\prod_{i=1}^n 2^{-\ve t|h_i|}~\C_{p~q~\gamma~\eta}.
\end{array}
\eeq
The estimate in 
(\ref{U_h(t)}) with $v=0$ implies  (\ref{Ortho Result}) for every $0<\alpha<n$.

Let  $q\in\Z,~\left({q-2\over p}\right)\left({q-r\over q}\right)\ge1$. 
By carrying out the same estimate in (\ref{U_h bound}) with $\Delta_\t\J_{\alpha~z} f$ replaced by $\Delta_\t\I_\alpha f$, we find 
\bel{Ortho Result each}
\begin{array}{lr}\ds
\int_{\R^n} \Big(\Delta_\t \I_\alpha f\Big)(x) \Big(\Delta_{\t-\h}\I_\alpha f\Big)^{q-1}(x)d\mu(x)
\\\\ \ds
~\leq~\C_{p~q~\gamma~\eta}\prod_{i=1}^n 2^{-\ve|h_i|}\left\| f\right\|_{\L^p(\mu)}^{q-r}\int_{\R^n} \Big(\vartheta_{\t-\h}(x)\Big)^{\left({q-2\over p}\right)\left({q-r\over q}\right)} \Big(\M_\beta f\Big)^r(x)d\mu (x)
\end{array}
\eeq
for  some $\ve=\ve(p,q,\gamma,\eta)>0$ and every  $\t$ and $\h$. Note that (\ref{Ortho Result}) can be deduced from (\ref{Ortho Result each}).

\v

{\bf ( 3 )}  It is suffice to prove  (\ref{Ortho Result each}) for $\t=\h$. ( i.e: $\t-\h=o$ )

Write 
\bel{Dila s}
\s x~\doteq~\left(2^{-s_1}x_1,~2^{-s_2}x_2,~\ldots,~2^{-s_n}x_n\right),\qquad
\s^{-1} x~\doteq~\left(2^{s_1}x_1,~2^{s_2}x_2,~\ldots,~2^{s_n}x_n\right)
\eeq
for  $s_i\in\Z,~i=1,2,\ldots,n $. 

Define
\bel{V_s mu_s}
d\mu_\s(x)~\doteq~\omega(\s x)dx,\qquad\V_\s(x,y)~\doteq~\inf_\delta~\Bigg\{ \mu_\s\left\{\Q( x,\delta)\right\}~\colon~ y\in\Q( x,\delta)\Bigg\}.
\eeq
We can  easily verify 
\bel{dila derivatives}
\begin{array}{cc}\ds
d\mu(\s x)~=~\prod_{i=1}^n 2^{-s_i}~\omega(\s x)dx~=~\prod_{i=1}^n 2^{-s_i}~ d\mu_\s(x),
\\\\ \ds
\V(\s x,\s y)~=~\prod_{i=1}^n 2^{-s_i}~ \V_\s(x,y).
\end{array}
\eeq
Note that $\mu_\s$ is {\it rectangle~doubling} with the same exponents $\gamma$ and $\eta(\gamma)$ in (\ref{rectangle doubling}). 

Recall $\Gamma_\t(x)$ defined in (\ref{Cone}). 
Let $f_\s(x)\doteq f(\s x)$. From (\ref{V_s mu_s})-(\ref{dila derivatives}), we have
\bel{Partial Dila}
\begin{array}{lr}\ds
\Big(\Delta_\t \I_\alpha f\Big)(\s x)~=~ \int_{\Gamma_{\t}(\s x)} f( y)\left({1\over \V(\s x, y)}\right)^{1-{\alpha\over n}} d\mu(y)
\\\\ \ds~~~~~~~~~~~~~~~~~~~
~=~ \int_{\Gamma_{\t}(\s x)} f( \s y)\left({1\over \V(\s x, \s y)}\right)^{1-{\alpha\over n}} d\mu(\s y)\qquad ( y\mt\s y )
\\\\ \ds~~~~~~~~~~~~~~~~~~~
~=~\prod_{i=1}^n 2^{-s_i\left({\alpha\over n}\right)} \int_{\Gamma_{\t-\s}(x)} f_\s( y)\left({1\over \V_\s( x, y)}\right)^{1-{\alpha\over n}} d\mu_\s(y)
\\\\ \ds~~~~~~~~~~~~~~~~~~~
~\leq~\C_{p~q~\gamma~\eta} \prod_{i=1}^n 2^{-s_i\left({\alpha\over n}\right)}\Big(\vartheta_{\t-\s}^\s(x)\Big)^{{1\over p}\left(1-{r\over q}\right)}\Big(\M_\beta^\s f_\s\Big)^{r\over q}(x)\left\| f_\s\right\|_{\L^p\left(\mu_\s\right)}^{1-{r\over q}}
\end{array}
\eeq
where $\vartheta_{\t-\s}^\s(x)$ and $\M_\beta^\s f_\s$ are defined as (\ref{theta_t}) and (\ref{M_alpha}) with $\mu$ and $f$ replaced by $\mu_\s$ and $f_\s$ respectively.

In particular, we have
\bel{theta invariant}
\begin{array}{lr}\ds
\vartheta_{\t-\s}^\s( x)~\doteq~      { \ds  \int_{\Gamma_{\t-\s}(x)} \Big(f_\s(y)\Big)^p d\mu_\s(y)\over \ds \int_{\R^n} \Big(f_\s(y)\Big)^pd\mu_\s(y)}
~=~{ \ds \prod_{i=1}^n 2^{s_i} \int_{\Gamma_{\t-\s}(x)} \Big(f(\s y)\Big)^p d\mu(\s y)\over \ds \prod_{i=1}^n 2^{s_i} \int_{\R^n} \Big(f( \s y)\Big)^pd\mu( \s y)}
\\\\ \ds~~~~~~~~~~~
~=~{ \ds  \int_{\Gamma_{\t-\s}(x)} \Big(f(\s y)\Big)^p d\mu(\s y)\over \ds  \int_{\R^n} \Big(f( \s y)\Big)^pd\mu( \s y)}~=~{ \ds  \int_{\Gamma_{\t}(\s x)} \Big(f( y)\Big)^p d\mu( y)\over \ds  \int_{\R^n} \Big(f(  y)\Big)^pd\mu(  y)} \qquad ( \s y\mt y )
\\\\ \ds~~~~~~~~~~~
~\doteq~\vartheta_\t(\s x)
\end{array}
\eeq
and
\bel{M Dila}
\begin{array}{lr}\ds
\Big(\M_\beta^\s f_\s\Big)(x)~=~\sup_\delta  \left\{\int_{\Q(x,\delta)}d\mu_\s(y)\right\}^{{\beta\over n}-1}\int_{\Q(x,\delta)} f_\s(y)d\mu_\s(y)
\\\\ \ds~~~~~~~~~~~~~~~~
~=~\sup_\delta  \left\{\prod_{i=1}^n 2^{s_i} \int_{\Q( x,\delta)}d\mu(\s y)\right\}^{{\beta\over n}-1}\prod_{i=1}^n 2^{s_i}\int_{\Q(x,\delta)} f(\s y)d\mu(\s y)
\\\\ \ds~~~~~~~~~~~~~~~~
~=~\prod_{i=1}^n 2^{s_i\left({\beta\over n}\right)}  \sup_\delta  \left\{ \int_{\Q(\s  x,\delta)}d\mu( y)\right\}^{{\beta\over n}-1}\int_{\Q(\s x,\delta)} f( y)d\mu( y)\qquad ( \s y\mt y )
\\\\ \ds~~~~~~~~~~~~~~~~
~=~\prod_{i=1}^n 2^{s_i\left({\beta\over n}\right)}\Big(\M_\beta f\Big)(\s x).
\end{array}
\eeq
By using (\ref{dila derivatives})-(\ref{M Dila}) with $\s=\t-\h$,   we have
\bel{Dila Est sum t}
\begin{array}{lr}\ds
\int_{\R^n}\Big(\Delta_\t \I_\alpha f\Big)(x) \Big(\Delta_{\t-\h}\I_\alpha f\Big)^{q-1}(x)d\mu(x)
\\\\ \ds
~=~ \int_{\R^n} \Big(\Delta_{\t} \I_\alpha f\Big)(\t-\h~ x) \Big(\Delta_{\t-\h}\I_\alpha f\Big)^{q-1}(\t-\h~ x)d\mu(\t-\h~ x)
\qquad(~x~\mt~\t-\h~ x~)
\\\\ \ds
~\leq~\C_{p~q~\gamma~\eta}\prod_{i=1}^n 2^{-(t_i-h_i)q\left({\alpha\over n} +{1\over q}\right)}
\prod_{i=1}^n 2^{-\ve|h_i|}\left\| f_{\t-\h}\right\|_{\L^p\left(\mu_{\t-\h}\right)}^{q-r}
 \int_{\R^n}\vartheta_{o}^{\t-\h}\left(x\right)^{\left({q-2\over p}\right)\left({q-r\over q}\right)}\Big(\M_{\beta}^{\t-\h} f_{\t-\h}\Big)^r(x)d\mu_{\t-\h} (x)
\\\\ \ds
~=~\C_{p~q~\gamma~\eta}
\prod_{i=1}^n 2^{(t_i-h_i)\left[\left({q-r\over p}\right)+r\left({\beta\over n}\right)+1-q\left({\alpha\over n} +{1\over q}\right)\right] }
\\\\ \ds~~~~~~~
\prod_{i=1}^n 2^{-\ve|h_i|}\left\| f\right\|_{\L^p\left(\mu\right)}^{q-r}\int_{\R^n}\vartheta_{\t-\h}\left(\t-\h ~x\right)^{\left({q-2\over p}\right)\left({q-r\over q}\right)}\Big(\M_{\beta} f\Big)^r(\t-\h~x)d\mu (\t-\h~x)
\\\\ \ds
~=~\C_{p~q~\gamma~\eta}\prod_{i=1}^n 2^{-\ve|h_i|}\left\| f\right\|_{\L^p\left(\mu\right)}^{q-r}\int_{\R^n}\vartheta_{\t-\h}\left(\t-\h ~x\right)^{\left({q-2\over p}\right)\left({q-r\over q}\right)}\Big(\M_{\beta} f\Big)^r(\t-\h~x)d\mu (\t-\h~x)
\end{array}
\eeq
where ${\alpha\over n}={1\over p}-{1\over q},~{\beta\over n}={1\over p}-{1\over r}$ implies $\left({q-r\over p}\right)+r\left({\beta\over n}\right)+1-q\left({\alpha\over n} +{1\over q}\right)=0$.

\section{Proof of the almost orthogonality}
\setcounter{equation}{0}
Let ${\alpha\over n}={1\over p}-{1\over q},1<p<q<\infty$ and  ${\beta\over n}={1\over p}-{1\over r},~p<r<q$. From {\bf (~1~)}-{\bf (~3~)} in the previous section, we aim to show
\bel{Est Orthogonality}
\begin{array}{lr}\ds
\int_{\R^n} \Big(\Delta_\t \I_\alpha f\Big)(x) \Big(\Delta_{o}\I_\alpha f\Big)^{q-1}(x)d\mu(x)
\\\\ \ds
~\leq~\C_{p~q~\gamma~\eta}\prod_{i=1}^n 2^{-\ve|t_i|}\left\| f\right\|_{\L^p(\mu)}^{q-r}\int_{\R^n} \Big(\vartheta_{o}(x)\Big)^{\left({q-2\over p}\right)\left({q-r\over q}\right)} \Big(\M_\beta f\Big)^r(x)d\mu (x)
\end{array}
\eeq
for some $\ve=\ve(p,q,\gamma,\eta)>0$ and every  $\t$, provided that
\bel{Constraints}
q\in\Z,\qquad\left({q-2\over p}\right)\left({q-r\over q}\right)\ge1,\qquad \eta-\gamma\left(1-{\alpha\over n}\right)>0.
\eeq 
From (\ref{Partial})-(\ref{Cone}) and (\ref{Gamma_t^j}), we write
\bel{Est explicite}
\begin{array}{lr}\ds
\int_{\R^n} \Big(\Delta_\t \I_\alpha f\Big)(x) \Big(\Delta_{o}\I_\alpha f\Big)^{q-1}(x)d\mu(x)~=~
\\\\ \ds
\int_{\R^n}\sum_{ \jmath,\ell_1,\ell_2,\ldots,\ell_{q-1}\in\Z}
\left\{\int_{\Gamma_{\t}^\jmath(x)} f( y)\left({1\over \V( x, y)}\right)^{1-{\alpha\over n}} d\mu(y)\right\} \prod_{m=1}^{q-1}\left\{\int_{\Gamma_{o}^{\ell_m}(x)} f( z^m)\left({1\over \V( x, z^m)}\right)^{1-{\alpha\over n}} d\mu(z^m)\right\} d\mu(x).
\end{array}
\eeq

Let
\bel{t_imath ell min}
\begin{array}{cc}
t_\imath~=~\max\{t_i\colon i=1,2,\ldots,n\},
\qquad
\ell_\nu~=~\min\{ \ell_m\colon m=1,2,\ldots,q-1\}.
\end{array}
\eeq
We develop a $3$-fold estimate by splitting the sum in (\ref{Est explicite}) into three groups:
\bel{Splitting G}
\sum_{ \jmath,\ell_1,\ell_2,\ldots,\ell_{q-1}~\in~\Z}~=~\sum_{\G_1}+\sum_{\G_2}+\sum_{\G_3} ~~~~,
\eeq
\bel{groups}
\begin{array}{lr}\ds
\G_1~=~\left\{ \jmath, \ell_1,\ell_2,\ldots,\ell_{q-1}\in\Z\colon~\jmath-t_\imath~\ge~\ell_\nu-2\right\},
\qquad
\G_2~=~\left\{ \jmath, \ell_1,\ell_2,\ldots,\ell_{q-1}\in\Z\colon~\jmath~\leq~\ell_\nu\right\},
\\\\ \ds~~~~~~~~~~~~~~~~~~~~~
\G_3~=~\left\{ \jmath, \ell_1,\ell_2,\ldots,\ell_{q-1}\in\Z\colon~\jmath-t_\imath~<~\ell_\nu-2~<~\jmath-2\right\}.
\end{array}
\eeq
Recall (\ref{tau(t,x)}). We have
\bel{tau, formula m,t}
\begin{array}{lr}\ds
\vartheta_\t(x)~=~{\Big(\M_\beta f\Big)^p(x)\over \left\| f\right\|_{\L^p\left(\mu\right)}^p}\mu\left\{\bigotimes_{i=1}^n |x_i-y_i|\leq2^{\tau(\t,x)-t_i}\right\}^{p\over r}
\end{array}
\eeq
where $0\leq\vartheta_\t(x)\leq1$ is defined in (\ref{theta_t}).

Let
\bel{tau(x)} 
\tau(x)~=~\tau(\t,x)\qquad\hbox{for}\qquad t_1=t_2=\cdots=t_n=0.
\eeq
Denote $j$ and $l_m,~ m=1,2,\ldots,q-1$ implicitly by
\bel{j,l_m}
\jmath~=~\tau(\t,x)+j\qquad\hbox{and}\qquad \ell_m~=~\tau(x)+l_m,~~ m=1,2,\ldots, q-1.
\eeq
\v

{\bf Case 1:} Consider $\jmath-t_\imath\ge\ell_\nu-2$. 
Suppose $\tau(\t,x)-\tau(x)>(1-\sigma) t_\imath$ for some $\sigma>0$. We have
\bel{Compara Est1}
\begin{array}{lr}\ds
{\vartheta_o(x)\over\vartheta_\t(x)}~\leq~{\ds\mu\left\{\bigotimes_{i=1}^n |x_i-y_i|\leq2^{\tau(x)}\right\}^{p\over r}\over\ds\mu\left\{|x_\nu-y_\nu|\leq2^{\tau(\t,x)}\times\bigotimes_{i\neq\nu} 
|x_i-y_i|\leq2^{\tau(\t,x)-t_\imath}\right\}^{p\over r}}
\\\\ \ds~~~~~~~~
~=~{\ds\mu\left\{\bigotimes_{i=1}^n |x_i-y_i|\leq2^{\tau(x)}\right\}^{p\over r}\over\ds\mu\left\{|x_\nu-y_\nu|\leq2^{\tau(\t,x)}\times\bigotimes_{i\neq\nu} |x_i-y_i|\leq2^{\tau(x)}\right\}^{p\over r}}
{\ds\mu\left\{|x_\nu-y_\nu|\leq2^{\tau(\t,x)}\times\bigotimes_{i\neq\nu} |x_i-y_i|\leq2^{\tau(x)}\right\}^{p\over r}\over\ds \mu\left\{|x_\nu-y_\nu|\leq2^{\tau(\t,x)}\times\bigotimes_{i\neq\nu} |x_i-y_i|\leq2^{\tau(\t,x)-t_\imath}\right\}^{p\over r}}
\\\\ \ds~~~~~~~
~\leq~2^{-\eta(p/r)(1-\sigma)t_\imath}\prod_{i\neq\nu} 2^{\gamma(p/r)\sigma t_\imath}
~=~2^{-(p/r)\left(\eta(1-\sigma)-\gamma\sigma(n-1)\right)t_\imath}.
\end{array}
\eeq
Note that $\eta(1-\sigma)-\gamma\sigma(n-1)>0$ for $\sigma$  sufficiently small.

By using (\ref{Regularity Est Rectangle}), we have
\bel{Ortho Result H<}
\begin{array}{lr}\ds
\int_{\Gamma_\t^\jmath(x)}f(y)\left({1\over\V(x,y)}\right)^{1-{\alpha\over n}} d\mu(y)
\prod_{m=1}^{q-1}\int_{\Gamma_o^{\ell_m}(x)}f(z^m)\left({1\over\V(x,z^m)}\right)^{1-{\alpha\over n}} d\mu(z^m)
\\\\ \ds
~\leq~2^{q\gamma n} ~ 2^{-\eta\min\left\{\alpha-\beta, {n\over q} \right\}|j|} \prod_{m=1}^{q-1} 2^{-\eta\min\left\{\alpha-\beta, {n\over q} \right\}|l_m|} \Big(\vartheta_\t(x)\Big)^{{1\over p}\left(1-{r\over q}\right)}\Big(\vartheta_o(x)\Big)^{{1\over p}\left(1-{r\over q}\right)}
\\\\ \ds~~~~~~~~~~~~
\Big(\vartheta_o(x)\Big)^{\left({q-2\over p}\right)\left({q-r\over q}\right)}\Big(\M_\beta f\Big)^{r}(x)\left\| f\right\|_{\L^p\left(\mu\right)}^{q-r}
\\\\ \ds
~\leq~2^{q\gamma n}~2^{-\left(\eta(1-\sigma)- \gamma\sigma (n-1)\right)\left({1\over r}-{1\over q}\right)t_\imath}   2^{-\eta\min\left\{\alpha-\beta, {n\over q} \right\}|j|} \prod_{m=1}^{q-1} 2^{-\eta\min\left\{\alpha-\beta, {n\over q} \right\}|l_m|}  
\\\\ \ds~~~~~~~~~~~~
 \Big(\vartheta_o(x)\Big)^{\left({q-2\over p}\right)\left({q-r\over q}\right)}\Big(\M_\beta f\Big)^{r}(x)\left\| f\right\|_{\L^p\left(\mu\right)}^{q-r}\qquad\hbox{\small{by (\ref{Compara Est1})}}.
\end{array}
\eeq
On the other hand, suppose $\tau(\t,x)-\tau(x)\leq(1-\sigma) t_\imath$ for some $\sigma>0$.
We have
\bel{j>l est}
\jmath-t_\imath~=~\tau(\t,x)+j-t_\imath~\ge~\tau(x)+l_\nu-2~=~\ell_\nu-2
\eeq
which  implies
\bel{j,l differs >1}
\begin{array}{lr}
j-l_\nu~\ge~t_\imath-\Big(\tau(\t,x)-\tau(x)\Big)-2
~\ge~t_\imath-(1-\sigma) t_\imath-1~=~\sigma t_\imath-2.
\end{array}
\eeq
By using (\ref{Regularity Est Rectangle}), we have
\bel{Ortho Est case1}
\begin{array}{lr}\ds
\int_{\Gamma_\t^\jmath(x)}f(y)\left({1\over\V(x,y)}\right)^{1-{\alpha\over n}} d\mu(y)
\prod_{m=1}^{q-1}\int_{\Gamma_o^{\ell_m}(x)}f(z^m)\left({1\over\V(x,z^m)}\right)^{1-{\alpha\over n}} d\mu(z^m)
\\\\ \ds
~\leq~ 2^{q\gamma n}~ 2^{-\eta\min\left\{\alpha-\beta, {n\over q} \right\}|j|} \prod_{m=1}^{q-1} 2^{-\eta\min\left\{\alpha-\beta, {n\over q} \right\}|l_m|}  \Big(\vartheta_\t(x)\Big)^{{1\over p}\left(1-{r\over q}\right)}\Big(\vartheta_o(x)\Big)^{{1\over p}\left(1-{r\over q}\right)}
\\\\ \ds~~~~~~~~~~~~
\Big(\vartheta_o(x)\Big)^{\left({q-2\over p}\right)\left({q-r\over q}\right)}\Big(\M_\beta f\Big)^{r}(x)\left\| f\right\|_{\L^p\left(\mu\right)}^{q-r}
\\\\ \ds
~\leq~2^{q\gamma n}~ 2^{-{1\over 2}\eta\min\left\{\alpha-\beta, {n\over q} \right\} \left(j-l_\nu\right)} 2^{-{1\over 2}\eta\min\left\{\alpha-\beta, {n\over q} \right\}|j|} \prod_{m=1}^{q-1} 2^{-{1\over 2}\eta\min\left\{\alpha-\beta, {n\over q} \right\}|l_m|}  

\\\\ \ds~~~~~~~~~~~~
\Big(\vartheta_o(x)\Big)^{\left({q-2\over p}\right)\left({q-r\over q}\right)}\Big(\M_\beta f\Big)^{r}(x)\left\| f\right\|_{\L^p\left(\mu\right)}^{q-r}
\qquad( j-l_\nu\leq|j|+|l_\nu|)
\\\\ \ds
~\leq~\C_{p~q~\gamma~\eta}~  2^{-{1\over 2}\eta\min\left\{\alpha-\beta, {n\over q} \right\} \sigma t_\imath} 2^{-{1\over 2}\eta\min\left\{\alpha-\beta, {n\over q} \right\}|j|} \prod_{m=1}^{q-1} 2^{-{1\over 2}\eta\min\left\{\alpha-\beta, {n\over q} \right\}|l_m|}  

\\\\ \ds~~~~~~~~~~~~
\Big(\vartheta_o(x)\Big)^{\left({q-2\over p}\right)\left({q-r\over q}\right)}\Big(\M_\beta f\Big)^{r}(x)\left\| f\right\|_{\L^p\left(\mu\right)}^{q-r}\qquad\hbox{\small{by (\ref{j,l differs >1})}}.
\end{array}
\eeq

By putting together  (\ref{Ortho Result H<}) and (\ref{Ortho Est case1}), we find
\bel{Result Case1}
\begin{array}{lr}\ds
\int_{\R^n}\sum_{\G_1}
\left\{\int_{\Gamma_{\t}^\jmath(x)} f( y)\left({1\over \V( x, y)}\right)^{1-{\alpha\over n}} d\mu(y)\right\} \prod_{m=1}^{q-1}\left\{\int_{\Gamma_{o}^{\ell_m}(x)} f( z^m)\left({1\over \V( x, y)}\right)^{1-{\alpha\over n}} d\mu(z^m)\right\} d\mu(x)
\\\\ \ds
~\leq~\C_{p~q~\gamma~\eta} 2^{-\ve_1 t_\imath} \left\| f\right\|_{\L^p\left(\mu\right)}^{q-r}
\\\\ \ds~~~~~~~
 \int_{\R^n}\Bigg\{\sum_{j,l_1,l_2,\ldots,l_{q-1}} 2^{-{1\over 2}\eta\min\left\{\alpha-\beta, {n\over q} \right\}|j|} \prod_{m=1}^{q-1} 2^{-{1\over 2}\eta\min\left\{\alpha-\beta, {n\over q} \right\}|l_m|} \Bigg\}
\Big(\vartheta_o(x)\Big)^{\left({q-2\over p}\right)\left({q-r\over q}\right)}\Big(\M_\beta f\Big)^{r}(x) d\mu(x)
\\\\ \ds
~\leq~\C_{p~q~\gamma~\eta}\prod_{i=1}^n 2^{-(\ve_1/n)t_i}  \left\| f\right\|_{\L^p\left(\mu\right)}^{q-r}\int_{\R^n}\Big(\vartheta_o(x)\Big)^{\left({q-2\over p}\right)\left({q-r\over q}\right)}\Big(\M_\beta f\Big)^{r}(x) d\mu(x).
\end{array}
\eeq
The exponent $\ve_1$ in (\ref{Result Case1}) equals
\bel{ve_1}
\begin{array}{lr}\ds
  \Big(\eta(1-\sigma)- \gamma\sigma (n-1)\Big)\left({1\over r}-{1\over q}\right)~=~{1\over 2}\sigma\eta\min\left\{\alpha-\beta,~{n\over q}\right\}~>~0
\end{array}
\eeq
where $0<\sigma<1$ is chosen so that the  equality holds in (\ref{ve_1}).
\v

{\bf Case 2:} Consider $\jmath\leq\ell_\nu$. 
Suppose $\tau(\t,x)-\tau(x)<\sigma t_\imath$ for some $\sigma>0$. We have
\bel{Compara Est2}
\begin{array}{lr}\ds
{\vartheta_\t(x)\over\vartheta_o(x)}
~\leq~{\ds\mu\left\{|x_\imath-y_\imath|\leq2^{\tau(\t,x)-t_\imath}\times\bigotimes_{i\neq\imath} |x_i-y_i|\leq2^{\tau(\t,x)}\right\}^{p\over r}\over\ds\mu\left\{\bigotimes_{i=1}^n |x_i-y_i|\leq2^{\tau(x)}\right\}^{p\over r}}
\\\\ \ds~~~~~~~~~~
~=~{\ds\mu\left\{|x_\imath-y_\imath|\leq2^{\tau(\t,x)-t_\imath}\times\bigotimes_{i\neq\imath} |x_i-y_i|\leq2^{\tau(\t,x)}\right\}^{p\over r}\over\ds\mu\left\{|x_\imath-y_\imath|\leq2^{\tau(x)}\times\bigotimes_{i\neq\imath} |x_i-y_i|\leq2^{\tau(\t,x)}\right\}^{p\over r}}
\\\\ \ds~~~~~~~~~~~~
\times~{\ds\mu\left\{|x_\imath-y_\imath|\leq2^{\tau(x)}\times\bigotimes_{i\neq\imath} |x_i-y_i|\leq2^{\tau(\t,x)}\right\}^{p\over r}\over \ds\mu\left\{\bigotimes_{i=1}^n |x_i-y_i|\leq2^{\tau(x)}\right\}^{p\over r}}
\\\\ \ds~~~~~~~~~~
~\leq~2^{-\eta(p/r)(1-\sigma)t_\imath}\prod_{i\neq\imath} 2^{\gamma(p/r)\sigma t_\imath}
~=~2^{-(p/r)\left(\eta(1-\sigma)-\gamma\sigma(n-1)\right)t_\imath}.
\end{array}
\eeq
By using (\ref{Regularity Est Rectangle}), we have
\bel{Ortho Result H>}
\begin{array}{lr}\ds
\int_{\Gamma_\t^\jmath(x)}f(y)\left({1\over\V(x,y)}\right)^{1-{\alpha\over n}} d\mu(y)
\prod_{m=1}^{q-1}\int_{\Gamma_o^{\ell_m}(x)}f(z^m)\left({1\over\V(x,z^m)}\right)^{1-{\alpha\over n}} d\mu(z^m)
\\\\ \ds
~\leq~2^{q\gamma n} ~ 2^{-\eta\min\left\{\alpha-\beta, {n\over q} \right\}|j|} \prod_{m=1}^{q-1} 2^{-\eta\min\left\{\alpha-\beta, {n\over q} \right\}|l_m|} \Big(\vartheta_\t(x)\Big)^{{1\over p}\left(1-{r\over q}\right)}\Big(\vartheta_o(x)\Big)^{{1\over p}\left(1-{r\over q}\right)}
\\\\ \ds~~~~~~~~~~~~
\Big(\vartheta_o(x)\Big)^{\left({q-2\over p}\right)\left({q-r\over q}\right)}\Big(\M_\beta f\Big)^{r}(x)\left\| f\right\|_{\L^p\left(\mu\right)}^{q-r}
\\\\ \ds
~\leq~2^{q\gamma n}~2^{-\left(\eta(1-\sigma)- \gamma\sigma (n-1)\right)\left({1\over r}-{1\over q}\right)t_\imath}   2^{-\eta\min\left\{\alpha-\beta, {n\over q} \right\}|j|} \prod_{m=1}^{q-1} 2^{-\eta\min\left\{\alpha-\beta, {n\over q} \right\}|l_m|}  
\\\\ \ds~~~~~~~~~~~~
 \Big(\vartheta_o(x)\Big)^{\left({q-2\over p}\right)\left({q-r\over q}\right)}\Big(\M_\beta f\Big)^{r}(x)\left\| f\right\|_{\L^p\left(\mu\right)}^{q-r}\qquad\hbox{\small{by (\ref{Compara Est2})}}.
\end{array}
\eeq
On the other hand, suppose $\tau(\t,x)-\tau(x)\ge\sigma t_\imath$. We have
\bel{j<l est}
\jmath~=~\tau(\t,x)+j~\leq~\tau(x)+l_\nu~=~\ell_\nu
\eeq
which  implies
\bel{j,l differs >2}
l_\nu-j~\ge~\tau(\t,x)-\tau(x)~\ge~\sigma t_\imath.
\eeq
By using (\ref{Regularity Est Rectangle}), we have 
\bel{Ortho Est case2}
\begin{array}{lr}\ds
\int_{\Gamma_\t^\jmath(x)}f(y)\left({1\over\V(x,y)}\right)^{1-{\alpha\over n}} d\mu(y)
\prod_{m=1}^{q-1}\int_{\Gamma_o^{\ell_m}(x)}f(z^m)\left({1\over\V(x,z^m)}\right)^{1-{\alpha\over n}} d\mu(z^m)
\\\\ \ds
~\leq~ 2^{q\gamma n}~ 2^{-\eta\min\left\{\alpha-\beta, {n\over q} \right\}|j|} \prod_{m=1}^{q-1} 2^{-\eta\min\left\{\alpha-\beta, {n\over q} \right\}|l_m|}  \Big(\vartheta_\t(x)\Big)^{{1\over p}\left(1-{r\over q}\right)}\Big(\vartheta_o(x)\Big)^{{1\over p}\left(1-{r\over q}\right)}
\\\\ \ds~~~~~~~~~~~~
\Big(\vartheta_o(x)\Big)^{\left({q-2\over p}\right)\left({q-r\over q}\right)}\Big(\M_\beta f\Big)^{r}(x)\left\| f\right\|_{\L^p\left(\mu\right)}^{q-r}
\\\\ \ds
~\leq~2^{q\gamma n}~ 2^{-{1\over 2}\eta\min\left\{\alpha-\beta, {n\over q} \right\} \left(l_\nu-j\right)} 2^{-{1\over 2}\eta\min\left\{\alpha-\beta, {n\over q} \right\}|j|} \prod_{m=1}^{q-1} 2^{-{1\over 2}\eta\min\left\{\alpha-\beta, {n\over q} \right\}|l_m|}  

\\\\ \ds~~~~~~~~~~~~
\Big(\vartheta_o(x)\Big)^{\left({q-2\over p}\right)\left({q-r\over q}\right)}\Big(\M_\beta f\Big)^{r}(x)\left\| f\right\|_{\L^p\left(\mu\right)}^{q-r}
\qquad( l_\nu-j\leq|j|+|l_\nu|)
\\\\ \ds
~\leq~2^{q\gamma n}~  2^{-{1\over 2}\eta\min\left\{\alpha-\beta, {n\over q} \right\} \sigma t_\imath} 2^{-{1\over 2}\eta\min\left\{\alpha-\beta, {n\over q} \right\}|j|} \prod_{m=1}^{q-1} 2^{-{1\over 2}\eta\min\left\{\alpha-\beta, {n\over q} \right\}|l_m|}  

\\\\ \ds~~~~~~~~~~~~
\Big(\vartheta_o(x)\Big)^{\left({q-2\over p}\right)\left({q-r\over q}\right)}\Big(\M_\beta f\Big)^{r}(x)\left\| f\right\|_{\L^p\left(\mu\right)}^{q-r}\qquad\hbox{\small{by (\ref{j,l differs >2})}}.
\end{array}
\eeq

By putting together (\ref{Ortho Result H>}) and (\ref{Ortho Est case2}), we find
\bel{Result Case2}
\begin{array}{lr}\ds
\int_{\R^n}\sum_{\G_1}
\left\{\int_{\Gamma_{\t}^\jmath(x)} f( y)\left({1\over \V( x, y)}\right)^{1-{\alpha\over n}} d\mu(y)\right\} \prod_{m=1}^{q-1}\left\{\int_{\Gamma_{o}^{\ell_m}(x)} f( z^m)\left({1\over \V( x, y)}\right)^{1-{\alpha\over n}} d\mu(z^m)\right\} d\mu(x)
\\\\ \ds
~\leq~\C_{p~q~\gamma~\eta} 2^{-\ve_2 t_\imath} \left\| f\right\|_{\L^p\left(\mu\right)}^{q-r}
\\\\ \ds~~~~~~~
 \int_{\R^n}\Bigg\{\sum_{j,l_1,l_2,\ldots,l_{q-1}} 2^{-{1\over 2}\eta\min\left\{\alpha-\beta, {n\over q} \right\}|j|} \prod_{m=1}^{q-1} 2^{-{1\over 2}\eta\min\left\{\alpha-\beta, {n\over q} \right\}|l_m|} \Bigg\}
\Big(\vartheta_o(x)\Big)^{\left({q-2\over p}\right)\left({q-r\over q}\right)}\Big(\M_\beta f\Big)^{r}(x) d\mu(x)
\\\\ \ds
~\leq~\C_{p~q~\gamma~\eta}\prod_{i=1}^n 2^{-(\ve_2/n)t_i}  \left\| f\right\|_{\L^p\left(\mu\right)}^{q-r}\int_{\R^n}\Big(\vartheta_o(x)\Big)^{\left({q-2\over p}\right)\left({q-r\over q}\right)}\Big(\M_\beta f\Big)^{r}(x) d\mu(x)
\end{array}
\eeq
where 
\bel{ve_2}
\begin{array}{lr}\ds
\ve_2~=~  \Big(\eta(1-\sigma)- \gamma\sigma (n-1)\Big)\left({1\over r}-{1\over q}\right)~=~{1\over 2}\sigma\eta\min\left\{\alpha-\beta,~{n\over q}\right\}~=~\ve_1~>~0.
\end{array}
\eeq
\v

{\bf Case 3:} Consider $\jmath-t_\imath<\ell_\nu-2<\jmath-2$.
Let $\{1,2,\ldots,n\}=\mathcal{U}\cup\mathcal{V}$ such that
\bel{UV}
\begin{array}{lr}\ds
\mathcal{U}~=~\left\{i\in\{1,2,\ldots,n\}\colon~\jmath+1-t_i\leq \ell_\nu-2\right\},
\qquad
\mathcal{V}~=~\left\{i\in\{1,2,\ldots,n\}\colon~\jmath+1-t_i> \ell_\nu-2\right\}.
\end{array}
\eeq
Note that $\imath\in\mathcal{U}$ because 
$t_\imath=\max \left\{t_i,i=1,2,\ldots,n\right\}$ in (\ref{t_imath ell min}).

Recall (\ref{Gamma_t^j}). We further write
\bel{rectangle t,j}
\begin{array}{lr}\ds
\Gamma_\t^{\jmath}(x)~=~\bigotimes_{i=1}^n\Gamma_{\t}^{\jmath~i}(x_i),
,\qquad
\Gamma_{\t}^{\jmath~i}(x_i)~=~\left\{y_i\in\R\colon~2^{\jmath-t_i}\leq |x_i-y_i|< 2^{\jmath+1-t_i}\right\}
\end{array}
\eeq
and their dyadic variants
\bel{rectangle* t,j}
\begin{array}{lr}\ds
{^\ast}\Gamma_{\t}^{\jmath}(x)~=~\bigotimes_{i=1}^n {^\ast}\Gamma_{\t}^{\jmath~i}(x_i),
,\qquad
{^\ast}\Gamma_{\t}^{\jmath~i}(x_i)~=~
\left\{x_i\in\R\colon~2^{\jmath-3-t_i}\leq |x_i-y_i|< 2^{\jmath+3-t_i}\right\}.
\end{array}
\eeq
From direct computation, we have
\bel{Ortho Expan}
\begin{array}{lr}\ds
~~~~~~~\int_{\R^n}\left\{\int_{\Gamma_\t^\jmath(x)}f(y)\left({1\over\V(x,y)}\right)^{1-{\alpha\over n}} d\mu(y)\right\}
\prod_{m=1}^{q-1}\left\{\int_{\Gamma_o^{\ell_m}(x)}f(z^m)\left({1\over\V(x,z^m)}\right)^{1-{\alpha\over n}} d\mu(z^m)\right\}d\mu(x)
\\\\ \ds
~=~\idotsint_{\R^n\times\cdots\times\R^n} f(y)\prod_{m=1}^{q-1}f(z^m)\left\{\int_{\Gamma_\t^\jmath(y)\cap\left\{\Cap_{m=1}^{q-1}\Gamma_o^{\ell_m}(z^m)\right\}}\left({1\over\V(x,y)}\right)^{1-{\alpha\over n}}\prod_{m=1}^{q-1}\left({1\over\V(x,z^m)}\right)^{1-{\alpha\over n}} d\mu(x)\right\} 
\\\\ \ds~~~~~~~~~~~~~~~~~~~~~~~~~~~~~~~~~~~~~~~~~~~~~~~~~~~~~~~~~~~~~~~~~~~~~~~~~~~~~~~~~~~~~~~~~~~~~~~~~~~~~
d\mu(y)\prod_{m=1}^{q-1}d\mu(z^m).
\end{array}
\eeq
Essentially,  we consider
\bel{Intersection} 
\Gamma^{\jmath}_\t(y)\cap\left\{ \Cap_{m=1}^{q-1} \Gamma^{\ell_m}_o(z^m)\right\}~\neq~\emptyset
\eeq
for $y,z^1,z^2\ldots,z^{q-1}\in\R^n$. 

Let $\r=\r(\t,\jmath-\ell_\nu)$ denote for the $n$-tuple $\left(2^{-r_1},~2^{-r_2},~\ldots,~2^{-r_n}\right)$ of which
\bel{r}
 \left\{\begin{array}{lr}\ds
r_i=\jmath-\ell_\nu+2,\qquad i\in\mathcal{U},
\\\\ \ds
r_i=t_i,\qquad i\in\mathcal{V}.
\end{array}\right.
\eeq
\begin{prop} There exists a cube
\bel{Inclusion}
\mathfrak{Q}~~~~\subset~~~~{^\ast}\Gamma_{\r}^{\jmath}(y_i)\cap\Bigg\{\Cap_{m=1}^{q-1}{^\ast}\Gamma_o^{\ell_m}(z^m)\Bigg\}
\eeq
such that
\bel{S norm}
\mu\left\{\bigotimes_{i\in\mathcal{U}}\Gamma_{o}^{\ell_\nu~i}(y_i)\times\bigotimes_{i\in\mathcal{V}}\Gamma_o^{\ell_\nu~i}(z^\nu_i)\right\}~\leq~\C_\gamma\mu\left\{\mathfrak{Q}\right\}
\eeq
whenever (\ref{Intersection}) holds for $y,z^1,z^2\ldots,z^{q-1}\in\R^n$.
\end{prop}
{\bf Proof:}  Consider $i\in\mathcal{U}$. 
From (\ref{UV}), we have $\jmath-t_i+1\leq\ell_\nu-2$. By  (\ref{Intersection}), there is an $\Hat{x}_i\in \Gamma^{\jmath~i}_\t(y_i)\cap\Big(\Cap_{m=1}^{q-1}\Gamma_o^{\ell_m~i}(z^m_i)\Big)$ such that 
\bel{norm compara}
\left|y_i-\Hat{x}_i\right|<2^{\jmath-t_i+1}\leq2^{\ell_\nu-2},\qquad 2^{\ell_m}\leq\left|\Hat{x}_i-z^m_i\right|<2^{\ell_m+1}.
\eeq 
By using the triangle inequality and (\ref{norm compara}), we have
\bel{norm Compara U centers}
\begin{array}{lr}\ds
2^{\ell_m-1}~<~2^{\ell_m}-2^{\ell_\nu-2}~<~\left|\Hat{x}_i-z^m_i\right|-\left|y_i-\Hat{x}_i\right|~\leq~\left|y_i-z^m_i\right|,
\\\\ \ds
\left|y_i-z^m_i\right|~\leq~\left|y_i-\Hat{x}_i\right|+\left|\Hat{x}_i-z^m_i\right|
~<~2^{\ell_\nu-2}+2^{\ell_m+1}~<~2^{\ell_m+2}.
\end{array}
\eeq
Let $x_i\in\Gamma_o^{\ell_\nu-3~i}(y_i)$ where $\left|x_i-y_i\right|<2^{\ell_\nu-2}$. By using the triangle inequality  and  (\ref{norm Compara U centers}), we  have
\bel{norm Compara U}
\begin{array}{lr}\ds
2^{\ell_m-3}~<~2^{\ell_m-1}-2^{\ell_\nu-2}~<~\left|y_i-z^m_i\right|-\left|x_i-y_i\right|~\leq~\left|x_i-z^m_i\right|,
\\\\ \ds
\left|x_i-z^m_i\right|~\leq~\left|x_i-y_i\right|+\left|y_i-z^m_i\right|~<~2^{\ell_\nu-2}+2^{\ell_m+2}~<~2^{\ell_m+3}.
\end{array}
\eeq
The estimates in (\ref{norm Compara U}) imply $x_i\in{^\ast}\Gamma_o^{\ell_m~i}(z^m_i)$ for every $m=1,2,\ldots,q-1$. Recall $\r$  defined in (\ref{r}).
We have $\Gamma_o^{\ell_\nu-2~i}(y_i)=\Gamma_\r^{\jmath~i}(y_i)\subset{^\ast}\Gamma_\r^{\jmath~i}(y_i)$ for $i\in\mathcal{U}$.

Consider $i\in\mathcal{V}$. From (\ref{UV}), we have $\jmath-t_i+1>\ell_\nu-2$. 
 Observe that $\Gamma_\t^{\jmath~i}(y_i)=\Gamma_\r^{\jmath~i}(y_i)$. Because
 of (\ref{Intersection}),  there is an $\Tilde{x}_i\in \Gamma^{\jmath~i}_\r(y_i)\cap\Big(\Cap_{m=1}^{q-1}\Gamma_o^{\ell_m~i}(z^m_i)\Big)$.

Let $\mathfrak{Q}_i\subset\Gamma_o^{\ell_\nu~i}(z^\nu_i)$  be a cube containing $\Tilde{x}_i$ whose side length equals $2^{\ell_\nu-3}$. It is clear that $\mathfrak{Q}_i$ intersects with $\Gamma^{\jmath~i}_\r(y_i)$ and every  $\Gamma_o^{\ell_m~i}(z_i^m),~m=1,2,\ldots,q-1$.
Recall (\ref{rectangle* t,j}). We must have
$\mathfrak{Q}_i\subset{^\ast}\Gamma_\r^{\jmath~i}(y_i)$, $\mathfrak{Q}_i\subset{^\ast}\Gamma_o^{\ell_m~i}(z_i^m), m=1,2,\ldots,q-1$
whereas $\ell_\nu-3<\jmath-t_i$ and $\ell_\nu-3<\ell_m$ for every $m=1,2,\ldots,q-1$.

Now, define $\mathfrak{Q}\doteq\bigotimes_{i\in\mathcal{U}} \Gamma_o^{\ell_\nu-3~i}(y_i)\times\bigotimes_{i\in\mathcal{V}}\mathfrak{Q}_i$  inside (\ref{Inclusion}). We write $2^5 \mathfrak{Q}$ for the cube having a same center of $\mathfrak{Q}$ but a $2^5$ times on its side length.   
As a geometric fact, we find $\bigotimes_{i\in\mathcal{U}}\Gamma_{o}^{\ell_\nu~i}(y_i)\times\bigotimes_{i\in\mathcal{V}}\Gamma_o^{\ell_\nu~i}(z^\nu_i)\subset 2^5 \mathfrak{Q}$.

The inequality holds in (\ref{S norm}) provided that $\mu$  is 
$rectangle~doubling$.\endproof

By using {\bf Proposition 7.1}, we have
\bel{Kernel Est}
\begin{array}{lr}\ds
\int_{\Gamma^{\jmath}_\t(y)\cap\left\{ \Cap_{m=1}^{q-1} \Gamma^{\ell_m}_o(z^m)\right\}}\left({1\over \V(x,y)}\right)^{1-{\alpha\over n}}\prod_{m=1}^{q-1}\left({1\over\V(x,z^m)}\right)^{1-{\alpha\over n}}d\mu(x)
\\\\ \ds
~\leq~\mu\left\{\bigotimes_{i=1}^n |x_i-y_i|\leq2^{\jmath-t_i}\right\}^{{\alpha\over n}-1}

\prod_{m=1}^{q-1}\mu\left\{\bigotimes_{i=1}^n |x_i-z_i^m|\leq2^{\ell_m}\right\}^{{\alpha\over n}-1}

\mu\left\{\bigotimes_{i\in\mathcal{U}} |x_i-y_i|\leq2^{\jmath-t_i}\times\bigotimes_{i\in\mathcal{V}} |x_i-z_i^\nu|\leq2^{\ell_\nu}     \right\}
\\\\ \ds

~\leq~\prod_{i\in\mathcal{U}}2^{\gamma\left(1-{\alpha\over n}\right)(t_i-\jmath+\ell_\nu)}\mu\left\{\bigotimes_{i\in\mathcal{U}} |x_i-y_i|\leq2^{\ell_\nu}\times\bigotimes_{i\in\mathcal{V}} |x_i-y_i|\leq2^{\jmath-t_i}\right\}^{{\alpha\over n}-1}
\\\\ \ds~~~~~~~
\prod_{i\in\mathcal{U}}2^{-\eta(t_i-\jmath+\ell_\nu)}\mu\left\{\bigotimes_{i\in\mathcal{U}} |x_i-y_i|\leq2^{\ell_\nu}\times\bigotimes_{i\in\mathcal{V}} |x_i-z_i^\nu|\leq2^{\ell_\nu}    \right\}
\prod_{m=1}^{q-1}\mu\left\{\bigotimes_{i=1}^n |x_i-z_i^m|\leq2^{\ell_m}\right\}^{{\alpha\over n}-1}
\\\\ \ds
~=~\prod_{i\in\mathcal{U}}2^{\left(\gamma\left(1-{\alpha\over n}\right)-\eta \right)(t_i-\jmath+\ell_\nu)}\mu\left\{\bigotimes_{i=1}^n |x_i-y_i|\leq2^{\jmath-r_i}\right\}^{{\alpha\over n}-1}
\prod_{m=1}^{q-1}\mu\left\{\bigotimes_{i=1}^n |x_i-z_i^m|\leq2^{\ell_m}\right\}^{{\alpha\over n}-1}
\\\\ \ds~~~~~~~
\mu\left\{\bigotimes_{i\in\mathcal{U}} |x_i-y_i|\leq2^{\ell_\nu}\times\bigotimes_{i\in\mathcal{V}} |x_i-z_i^\nu|\leq2^{\ell_\nu}      \right\}    
\qquad
\hbox{\small{by (\ref{r})}} 
\\\\ \ds
~\leq~ 2^{-\Big(\eta-\gamma\left(1-{\alpha\over n}\right)\Big)(t_\imath-\jmath+{\ell_\nu})} \C_\gamma \int_{{^\ast}\Gamma^{\jmath}_\r(y)\cap\left\{ \Cap_{m=1}^{q-1} {^\ast}\Gamma^{\ell_m}_o(z^m)\right\}}\left({1\over \V(x,y)}\right)^{1-{\alpha\over n}}\prod_{m=1}^{q-1}\left({1\over\V(x,z^m)}\right)^{1-{\alpha\over n}}d\mu(x)
\end{array}
\eeq
where the last inequality is obtained by using (\ref{Inclusion})-(\ref{S norm}). 

For $y\in{^\ast}\Gamma_\t^\jmath(x)$, we have $\V(x,y)$   bounded from both above and below as  shown in (\ref{V norm compara}) with two implied constants $2^{-3\eta n}$ and $2^{2\gamma n}$ added respectively. By carrying out (\ref{Region-I Est})-(\ref{Regularity Est II}) with $\Gamma^\jmath_\t(x)$ replaced by ${^\ast}\Gamma^\jmath_\t(x)$,
we find (\ref{Regularity Est Rectangle}) again except for  $ 2^{\gamma n}$  replaced by $\C_{p~q~\gamma~\eta}$.

Let $\r=\r(\t,\jmath-\ell_\nu)$ defined in (\ref{r}). Write
\bel{Hat j}
\jmath~=~j+\tau(\t,x)=\Hat{j}+\tau\Big(\r(\t,\jmath-\ell_\nu),x\Big)
\eeq 
for $\Hat{j}\in\Z-\tau\Big(\r(\t,\jmath-\ell_\nu),x\Big)$.
Moreover,  $\jmath-\ell_\nu=j-l_\nu+\Big(\tau(\t,x)-\tau(x)\Big)$. We denote
\bel{Theta}
\Theta(\t,x,j-l_\nu)~=~\tau(\t,x)-\tau\Big(\r(\t,\jmath-\ell_\nu),x\Big)
\eeq
which is a real number that depends only on $\t$, $x$ and $j-l_\nu$.

From (\ref{r}), we have $r_\imath=\min\{r_i\colon i=1,2,\ldots,n\}=\jmath-\ell_\nu-2$. Hence that $\jmath, \ell_1,\ell_2 \ldots, \ell_{q-1}$ belong to $\G_1$ defined in (\ref{Splitting G}) $w.r.t$ $\r$. 

By using  (\ref{Kernel Est}) and repeating all estimates  in {\bf Case 1} with  $\Gamma_\t^\jmath(x),~\Gamma_o^{\ell_m}(x),~m=1,2,\ldots,q-1$   replaced respectively by  ${^\ast}\Gamma_\r^\jmath(x),~{^\ast}\Gamma_o^{\ell_m}(x),~m=1,2,\ldots,q-1$,  we have
\bel{Ortho Expan Est Result}
\begin{array}{lr}\ds
\int_{\R^n}\left\{\int_{\Gamma_\t^\jmath(x)}f(y)\left({1\over\V(x,y)}\right)^{1-{\alpha\over n}} d\mu(y)\right\}
\prod_{m=1}^{q-1}\left\{\int_{\Gamma^{\ell_m}_o(x)}f(z^m)\left({1\over\V(x,z^m)}\right)^{1-{\alpha\over n}} d\mu(z^m)\right\}d\mu(x)
\\\\ \ds
~=~\idotsint_{\R^n\times\cdots\times\R^n} \left\{\int_{\Gamma_\t^\jmath(y)\cap\left\{\Cap_{m=1}^{q-1}\Gamma^{\ell_m}_o(z^m)\right\}}\left({1\over\V(x,y)}\right)^{1-{\alpha\over n}}\prod_{m=1}^{q-1}\left({1\over\V(x,z^m)}\right)^{1-{\alpha\over n}} d\mu(x)\right\} f(y)d\mu(y)\prod_{m=1}^{q-1}f(z^m)d\mu(z^m)
\\\\ \ds
~\leq~\C_\gamma~2^{-(\ve_3/2)(t_\imath+\ell_\nu-\jmath)}
\\\\ \ds
\idotsint_{\R^n\times\cdots\times\R^n} \left\{\int_{{^\ast}\Gamma_{\r}^\jmath(y)\cap\left\{\Cap_{m=1}^{q-1}{^\ast}\Gamma^{\ell_m}_o(z^m)\right\}}\left({1\over\V(x,y)}\right)^{1-{\alpha\over n}}\prod_{m=1}^{q-1}\left({1\over\V(x,z^m)}\right)^{1-{\alpha\over n}} d\mu(x)\right\} f(y)d\mu(y)\prod_{m=1}^{q-1}f(z^m)d\mu(z^m)
\\\\ \ds
~=~\C_\gamma~2^{-(\ve_3/2)(t_\imath+\ell_\nu-\jmath)}\int_{\R^n}\left\{\int_{{^\ast}\Gamma_{\r}^\jmath(x)}f(y)\left({1\over\V(x,y)}\right)^{1-{\alpha\over n}} d\mu(y)\right\}
\prod_{m=1}^{q-1}\left\{\int_{{^\ast}\Gamma^{\ell_m}_o(x)}f(z^m)\left({1\over\V(x,z^m)}\right)^{1-{\alpha\over n}} d\mu(z^m)\right\}d\mu(x)
\\\\ \ds
~\leq~\C_{p~q~\gamma~\eta}~ 2^{-(\ve_3/2)(t_\imath+\ell_\nu-\jmath)} 2^{-\ve_3|r_\imath|}\Bigg\{  2^{-{1\over 2}\eta\min\left\{\alpha-\beta, {n\over q} \right\}|\Hat{j}|} \prod_{m=1}^{q-1} 2^{-{1\over 2}\eta\min\left\{\alpha-\beta, {n\over q} \right\}|l_m|}  \Bigg\}
\\\\ \ds~~~~~~~~~~~~~~~~~~~~~~~~~~~~~~~~~~
\left\| f\right\|_{\L^p\left(\mu\right)}^{q-r}\int_{\R^n}\Big(\vartheta_o(x)\Big)^{\left({q-2\over p}\right)\left({q-r\over q}\right)}\Big(\M_\beta f\Big)^{r}(x)d\mu(x)
\\\\ \ds
~\leq~\C_{p~q~\gamma~\eta}~ 2^{-(\ve_3/2)(t_\imath+\ell_\nu-\jmath)} 2^{-\ve_3(\jmath-\ell_\nu)}\Bigg\{  2^{-{1\over 2}\eta\min\left\{\alpha-\beta, {n\over q} \right\}|\Hat{j}|} \prod_{m=1}^{q-1} 2^{-{1\over 2}\eta\min\left\{\alpha-\beta, {n\over q} \right\}|l_m|}  \Bigg\}
\\\\ \ds~~~~~~~~~~~~~~~~~~~~~~~~~~~~~~~~~~
\left\| f\right\|_{\L^p\left(\mu\right)}^{q-r}\int_{\R^n}\Big(\vartheta_o(x)\Big)^{\left({q-2\over p}\right)\left({q-r\over q}\right)}\Big(\M_\beta f\Big)^{r}(x)d\mu(x)
\\\\ \ds
~\leq~\C_{p~q~\gamma~\eta}~2^{-(\ve_3/2)t_\imath} 2^{-(\ve_3/2)(\jmath-\ell_\nu)}~  2^{-{1\over 2}\eta\min\left\{\alpha-\beta, {n\over q} \right\}\left|j+\tau(\t,x)-\tau\Big(\r(\t,\jmath-\ell_\nu),x\Big)\right|} \prod_{m=1}^{q-1} 2^{-{1\over 2}\eta\min\left\{\alpha-\beta, {n\over q} \right\}|l_m|}  
\\\\ \ds~~~~~~~~~~~~~~~~~~~~~~~~~~~~~~~~~~
\left\| f\right\|_{\L^p\left(\mu\right)}^{q-r}\int_{\R^n}\Big(\vartheta_o(x)\Big)^{\left({q-2\over p}\right)\left({q-r\over q}\right)}\Big(\M_\beta f\Big)^{r}(x)d\mu(x)
\\\\ \ds
~=~\C_{p~q~\gamma~\eta}~2^{-(\ve_3/2)t_\imath} 2^{-(\ve_3/2)\left|j-l_\nu+\Big(\tau(\t,x)-\tau(x)\Big)\right|}~2^{-{1\over 2}\eta\min\left\{\alpha-\beta, {n\over q} \right\}\left|j+\Theta(\t,x,j-l_\nu)\right|} \prod_{m=1}^{q-1} 2^{-{1\over 2}\eta\min\left\{\alpha-\beta, {n\over q} \right\}|l_m|}  
\\\\ \ds~~~~~~~~~~~~~~~~~~~~~~~~~~~~~~~~~~
\left\| f\right\|_{\L^p\left(\mu\right)}^{q-r}\int_{\R^n}\Big(\vartheta_o(x)\Big)^{\left({q-2\over p}\right)\left({q-r\over q}\right)}\Big(\M_\beta f\Big)^{r}(x)d\mu(x).
\end{array}
\eeq
The exponent $\ve_3$ in (\ref{Ortho Expan Est Result}) equals
\bel{ve_3}
\min\left\{\eta-\gamma\left(1-{\alpha\over n}\right),~ \Big(\eta(1-\sigma)- \gamma\sigma (n-1)\Big)\left({1\over r}-{1\over q}\right)\right\}={1\over 2}\sigma\eta
\min\left\{\alpha-\beta, {n\over q}\right\}~>~0
\eeq
where $0<\sigma<1$ is chosen to satisfy the equality in (\ref{ve_3}).

Lastly, by using (\ref{Ortho Expan Est Result})-(\ref{ve_3}), we have
\bel{Result Case3}
\begin{array}{lr}\ds
\int_{\R^n}\sum_{\G_3}
\left\{\int_{\Gamma_{\t}^\jmath(x)} f( y)\left({1\over \V( x, y)}\right)^{1-{\alpha\over n}} d\mu(y)\right\} \prod_{m=1}^{q-1}\left\{\int_{\Gamma_{o}^{\ell_m}(x)} f( z^m)\left({1\over \V( x, z^m)}\right)^{1-{\alpha\over n}} d\mu(z^m)\right\} d\mu(x)
\\\\ \ds
~=~\sum_{\G_3}\int_{\R^n}
\left\{\int_{\Gamma_{\t}^\jmath(x)} f( y)\left({1\over \V( x, y)}\right)^{1-{\alpha\over n}} d\mu(y)\right\} \prod_{m=1}^{q-1}\left\{\int_{\Gamma_{o}^{\ell_m}(x)} f( z^m)\left({1\over \V( x, z^m)}\right)^{1-{\alpha\over n}} d\mu(y)\right\} d\mu(x)
\\\\ \ds
~\leq~\C_{p~q~\gamma~\eta}~2^{-(\ve_3/2) t_\imath} \left\| f\right\|_{\L^p\left(\mu\right)}^{q-r}\sum_{\G_3}~\prod_{m=1}^{q-1} 2^{-{1\over 2}\eta\min\left\{\alpha-\beta, {n\over q} \right\}|l_m|}  
\\\\ \ds~~~~~~~
\int_{\R^n}2^{-(\ve_3/2)\left|j-l_\nu+\Big(\tau(\t,x)-\tau(x)\Big)\right|}2^{-{1\over 2}\eta\min\left\{\alpha-\beta, {n\over q} \right\}\left|j+\Theta(\t,x,j-l_\nu)\right|} \Big(\vartheta_o(x)\Big)^{\left({q-2\over p}\right)\left({q-r\over q}\right)}\Big(\M_\beta f\Big)^{r}(x)d\mu(x)
\\\\ \ds
~\leq~\C_{p~q~\gamma~\eta}~2^{-(\ve_3/2) t_\imath} \left\| f\right\|_{\L^p\left(\mu\right)}^{q-r}  \sum_{j,l_1,l_2,\ldots,l_{q-1}}    ~\prod_{m=1}^{q-1} 2^{-{1\over 2}\eta\min\left\{\alpha-\beta, {n\over q} \right\}|l_m|}    
\\\\ \ds~~~~~~~
\int_{\R^n}2^{-(\ve_3/2)\left|j-l_\nu+\Big(\tau(\t,x)-\tau(x)\Big)\right|}2^{-{1\over 2}\eta\min\left\{\alpha-\beta, {n\over q} \right\}\left|j+\Theta(\t,x,j-l_\nu)\right|} 
\Big(\vartheta_o(x)\Big)^{\left({q-2\over p}\right)\left({q-r\over q}\right)}\Big(\M_\beta f\Big)^{r}(x)d\mu(x)
\\\\ \ds
~\leq~\C_{p~q~\gamma~\eta}~2^{-(\ve_3/2) t_\imath} \left\| f\right\|_{\L^p\left(\mu\right)}^{q-r}\sum_{k, l_1,l_2,\ldots,l_{q-1}}  ~  \prod_{m=1}^{q-1} 2^{-{1\over 2}\eta\min\left\{\alpha-\beta, {n\over q} \right\}|l_m|}    
\\\\ \ds~~~~~~~
\int_{\R^n}2^{-(\ve_3/2)\left|k+\Big(\tau(\t,x)-\tau(x)\Big)\right|} \Big(\vartheta_o(x)\Big)^{\left({q-2\over p}\right)\left({q-r\over q}\right)}\Big(\M_\beta f\Big)^{r}(x)d\mu(x)
\qquad\hbox{\small{$(k=j-l_\nu)$}}
\\\\ \ds
~=~\C_{p~q~\gamma~\eta}~2^{-(\ve_3/2) t_\imath} \left\| f\right\|_{\L^p\left(\mu\right)}^{q-r}\sum_{ l_1,l_2,\ldots,l_{q-1}}    ~\prod_{m=1}^{q-1} 2^{-{1\over 2}\eta\min\left\{\alpha-\beta, {n\over q} \right\}|l_m|}   
\\\\ \ds~~~~~~~
\int_{\R^n}\left\{\sum_{k} 2^{-(\ve_3/2)\left|k+\Big(\tau(\t,x)-\tau(x)\Big)\right|}\right\}\Big(\vartheta_o(x)\Big)^{\left({q-2\over p}\right)\left({q-r\over q}\right)}\Big(\M_\beta f\Big)^{r}(x)d\mu(x)

\\\\ \ds
~\leq~\C_{p~q~\gamma~\eta} \prod_{i=1}^n 2^{-(\ve_3/2n) t_i} \left\| f\right\|_{\L^p\left(\mu\right)}^{q-r}
\int_{\R^n}
\Big(\vartheta_o(x)\Big)^{\left({q-2\over p}\right)\left({q-r\over q}\right)}\Big(\M_\beta f\Big)^{r}(x)d\mu(x).
\end{array}
\eeq

wangzipeng@westlake.edu.cn

\end{document}